\documentclass[12pt]{amsart} 
\usepackage[normalem]{ulem}  
\usepackage{lmodern}
\usepackage{amsmath}
\usepackage{amssymb}
\usepackage{dsfont}  % to do a doublestroke 1
\usepackage{hyperref}
\usepackage{graphicx}
\usepackage[all]{xy}
\usepackage{latexsym}
\usepackage{epstopdf}
\usepackage{pgf}
\usepackage{tikz}
\usepackage{tikz-cd}
\usepackage{pgf}
\usepackage{float}
\usepackage{subcaption} 
\usepackage{xcolor}
\usepackage{algorithm2e}
\usepackage{color}
\usepackage{ytableau}
\usepackage{adjustbox}
\usepackage{mathrsfs}
\usepackage{enumerate}   
\usepackage{mathbbol}
\usepackage{quiver}

\usetikzlibrary{matrix,arrows,backgrounds,shapes.misc,shapes.geometric,patterns,calc,positioning}
\usetikzlibrary{calc,shapes}
\usetikzlibrary{decorations.pathmorphing}
\usetikzlibrary{positioning,arrows}
\usetikzlibrary{decorations.markings}

\usepackage[left=3cm,right=3cm,top=2.8cm,bottom=2.8cm]{geometry}
\newtheorem{theorem}{Theorem}[section] 
\newtheorem{thm}{Theorem}
        
\newtheorem{lemma}[theorem]{Lemma}
\newtheorem{corollary}[theorem]{Corollary}

\newtheorem{proposition}[theorem]{Proposition}

\newtheorem{question}[thm]{Question}

\theoremstyle{remark}
\newtheorem{remark}[theorem]{Remark}

\theoremstyle{definition}
\newtheorem{definition}[theorem]{Definition} 
\newtheorem{example}[theorem]{Example} 

\numberwithin{equation}{section}

\newcommand{\rad}{\operatorname{rad}\nolimits}

\newcommand{\Hom}{\operatorname{Hom}\nolimits}

\newcommand{\Ext}{\operatorname{Ext}\nolimits}
\newcommand{\End}{\operatorname{End}\nolimits}

\newcommand{\opname}{\operatorname}

\newcommand{\dimv}{\underline{\mathrm{dim}}\hspace{1pt}}

\newcommand{\sqbinom}[2]{\begin{bmatrix} #1 \\ #2 \end{bmatrix}}

\begin{document}
\title{ Representations of bound quivers over $\mathbb{F}_1$ and Hall algebras}\thanks{This work was supported by the National Natural Science Foundation of China (Grant No. 12571040)}

\author{Changjian Fu}
\author{Xuelian Xiao}
\address{Department of Mathematics, Sichuan University, Chengdu, 610064 PR China}

\email{changjianfu@scu.edu.cn(Fu)}
\email{xiaoxuelian@stu.scu.edu.cn(Xiao)}

\begin{abstract}
    Representations of quivers over the virtual field $\mathbb{F}_1$ can be naturally extended to bound quivers with monomial and commutativity relations. In this paper, we introduce and study the Hall algebra $\mathcal{H}_{\mathbb{F}_1}(Q,I)$ for a bound quiver $(Q,I)$ equipped with such relations. As in the quiver setting, $\mathcal{H}_{\mathbb{F}_1}(Q,I)$ is a graded, connected, and cocommutative bialgebra, and the free subgroup generated by all indecomposable representations forms a Lie subalgebra. We study $\mathcal{H}_{\mathbb{F}_1}(Q,I)$ in detail for two special classes of bound quivers: Nakayama bound quivers and gentle one-cycle bound quivers. In both cases, we prove that for any triple $L, M, N$ of $\mathbb{F}_1$-representations over $(Q,I)$, the Hall number $g_{N,L}^M$ coincides with the evaluation of the corresponding Hall polynomial $\phi_{N,L}^M(t)$ at $t=1$. Consequently, there exists a canonical isomorphism between $\mathcal{H}_{\mathbb{F}_1}(Q,I)$ and the degenerate Hall algebra $\mathcal{H}_1(Q,I)$, justifying Tits's $\mathbb{F}_1$-philosophy in these settings.
\end{abstract}
\maketitle

\section{Introduction}

Let $\mathbb{F}_1$ be the virtual field, i.e., the ``field with one element''.
In \cite{Sz12}, Szczesny introduced quiver representations over $\mathbb{F}_1$ as a degenerate, combinatorial analogue of classical quiver representations over (finite) fields.  Over the past decade, the interplay between the combinatorics of $\mathbb{F}_1$-representations and their field-theoretic counterparts has been extensively explored (see, e.g., \cite{JKS23,JS23,JS24,Kleinau25}).

Let $Q$ be a finite quiver and $\operatorname{rep}_{\mathbb{F}_1}(Q)$ the category of representations of $Q$ over $\mathbb{F}_1$. The category $\operatorname{rep}_{\mathbb{F}_1}(Q)$ behaves in many respects like the category $\operatorname{rep}_k(Q)$ of representations of $Q$ over a field $k$, such as having a zero object, kernels, cokernels, and finite direct sums, as well as satisfying the Krull--Schmidt and Jordan--H\"{o}lder theorems. However, since $\operatorname{rep}_{\mathbb{F}_1}(Q)$ is not additive, it also exhibits significant differences from $\operatorname{rep}_k(Q)$. 
Indeed, as shown in \cite{FRY24, FYZ26b}, $\operatorname{rep}_{\mathbb{F}_1}(Q)$ does not necessarily have enough projectives, and its global dimension---defined via Yoneda extensions---can be $2$ (whereas the global dimension of $\operatorname{rep}_k(Q)$ is always at most $1$).

Nevertheless, the category $\operatorname{rep}_{\mathbb{F}_1}(Q)$ is a proto-exact category in the
sense of Dyckerhoff and Kapranov \cite{DK19}, where all monomorphisms and epimorphisms are
admissible. Consequently, one may define the Hall algebra $\mathcal{H}_{\mathbb{F}_1}(Q)$ of $\operatorname{rep}_{\mathbb{F}_1}(Q)$ in the same way as the Hall algebra $\mathcal{H}_{\mathbb{F}_q}(Q)$  of $\operatorname{rep}_{\mathbb{F}_q}(Q)$ over the finite field $\mathbb{F}_q$ (cf. \cite{Ringel90}). The Hall algebra $\mathcal{H}_{\mathbb{F}_1}(Q)$ behaves, in some sense, as a specialization at $q=1$ of $\mathcal{H}_{\mathbb{F}_q}(Q)$. In particular, when $Q$ has no loops,  Szczesny \cite{Sz12} proved that there is a Hopf algebra homomorphism $\rho: U(\mathfrak{n}_+)\to \mathcal{H}_{\mathbb{F}_1}(Q)$, where $U(\mathfrak{n}_+)$ is the enveloping algebra of the positive part $\mathfrak{n}_+$ of the Kac--Moody Lie algebra associated with $Q$. In contrast to the classical results for $\mathcal{H}_{\mathbb{F}_q}(Q)$ \cite{Ringel90,Green}, the morphism $\rho$ fails to be injective in general. On the other hand, the subspace $\mathfrak{n}_{\mathbb{F}_1}(Q) \subset \mathcal{H}_{\mathbb{F}_1}(Q)$ spanned by indecomposable representations naturally carries a Lie algebra structure, a property that does not hold for $\mathcal{H}_{\mathbb{F}_q}(Q)$ in general.

Motivated by these phenomena and aiming to move beyond the quiver setting, the goal of this paper is twofold. First, we observe that the definition of quiver representations over $\mathbb{F}_1$ extends naturally to bound quivers $(Q,I)$ with monomial and commutativity relations. The resulting category $\operatorname{rep}_{\mathbb{F}_1}(Q,I)$ remains a proto-exact category, allowing one to construct its associated Hall algebra $\mathcal{H}_{\mathbb{F}_1}(Q,I)$ and the Lie subalgebra $\mathfrak{n}_{\mathbb{F}_1}(Q,I)$ spanned by indecomposables. This broader framework yields a wealth of new Lie algebras. Second, it provides a suitable setting to justify Tits's philosophy on $\mathbb{F}_1$ \cite{Tits56}, namely that incidence geometries over a finite field $\mathbb{F}_q$ have combinatorial counterparts that can be interpreted as incidence geometries over $\mathbb{F}_1$.

Let $(Q,I)$ be a bound quiver with monomial and commutativity relations. For any finite field $\mathbb{F}_q$ and finite field extension $\mathbb{E}/\mathbb{F}_q$, there exist base change functors:
\begin{align*}
    -\otimes_{\mathbb{F}_1}\mathbb{F}_q:=-^{\mathbb{F}_q}:&\operatorname{rep}_{\mathbb{F}_1}(Q,I)\to \operatorname{rep}_{\mathbb{F}_q}(Q,I),\\
    -\otimes_{\mathbb{F}_q}\mathbb{E}:=-^{\mathbb{E}}:&\operatorname{rep}_{\mathbb{F}_q}(Q,I)\to \operatorname{rep}_{\mathbb{E}}(Q,I).
\end{align*}
For any $L,M,N\in \operatorname{rep}_{\mathbb{F}_1}(Q,I)$, denote
\[V_{\mathbb{F}_1}(L,M,N):=\{X\subseteq M\mid  X\cong L \text{ and } M/X\cong N \},\] 
\[V_{\mathbb{F}_q}(L^{\mathbb{F}_q},M^{\mathbb{F}_q},N^{\mathbb{F}_q}):=\{U\subseteq M^{\mathbb{F}_q}\mid U\cong L^{\mathbb{F}_q} \text{ and } M^{\mathbb{F}_q}/U\cong N^{\mathbb{F}_q} \}.\] 
According to Tits's philosophy, it is natural to ask whether $V_{\mathbb{F}_1}(L,M,N)$ is the combinatorial counterpart of the $\mathbb{F}_q$-variety $V_{\mathbb{F}_q}(L^{\mathbb{F}_q},M^{\mathbb{F}_q},N^{\mathbb{F}_q})$.
To formulate this question precisely, let us assume that $\operatorname{rep}_{\mathbb{F}_q}(Q,I)$ is representation-finite and admits Hall polynomials. In particular, there exists a polynomial $\phi_{N,L}^M(t)\in \mathbb{Z}[t]$ such that for any finite conservative extension $\mathbb{E}$ of $\mathbb{F}_q$, \[|V_\mathbb{E}(L^\mathbb{E},M^\mathbb{E},N^\mathbb{E})|=\phi_{N,L}^M(|\mathbb{E}|).\]
Furthermore, one can define the generic Hall algebra $\mathcal{H}_t(Q,I)$ and the degenerate Hall algebra $\mathcal{H}_1(Q,I)$ associated with $(Q,I)$.
The question can be phrased as follows:
\begin{question}\label{ques:tits-philosophy}
    \begin{itemize}
        \item[(1)] Do we have $|V_{\mathbb{F}_1}(L,M,N)|=\phi_{N,L}^M(1)$ for any triple $L,M,N\in \operatorname{rep}_{\mathbb{F}_1}(Q,I)$?
        \item[(2)] What is the relationship between the Hall algebra $\mathcal{H}_{\mathbb{F}_1}(Q,I)$ and the degenerate Hall algebra $\mathcal{H}_1(Q,I)$? What about their Lie subalgebras?
    \end{itemize}
\end{question}

In this paper, we address the above question for two classes of bound quivers: Nakayama bound quivers and gentle one-cycle bound quivers. The associated bound quiver algebras are string algebras, whose representation theories are well-established. In these cases, the base change functor $-\otimes_{\mathbb{F}_1}\mathbb{F}_q$ yields a bijection between the isomorphism classes of objects in $\operatorname{rep}_{\mathbb{F}_1}(Q,I)$ and $\operatorname{rep}_{\mathbb{F}_q}(Q,I)$. As our main results, we show that the equality in Question \ref{ques:tits-philosophy} (1) always holds, and we establish a canonical isomorphism between their Hall algebras. In particular, the Hall algebra $\mathcal{H}_{\mathbb{F}_1}(Q,I)$ is indeed the specialization at $t=1$ of the generic Hall algebra $\mathcal{H}_t(Q,I)$.
\begin{remark}
    For a general bound quiver $(Q,I)$ with monomial and commutativity relations, the category $\operatorname{rep}_{\mathbb{F}_q}(Q,I)$ may have  more indecomposables than  $\operatorname{rep}_{\mathbb{F}_1}(Q,I)$. Consequently, their corresponding Hall algebras need not be isomorphic. Nevertheless, whenever the degenerate Hall algebra $\mathcal{H}_1(Q,I)$ is well-defined, one might still expect the existence of a surjective homomorphism from $\mathcal{H}_1(Q,I)$ to $\mathcal{H}_{\mathbb{F}_1}(Q,I)$.
\end{remark}

The paper is organized as follows. In Section \ref{s:preliminary}, we recall the basics of Hall numbers and degenerate Hall algebras for representation-finite algebras over finite fields. We then provide definitions and foundational results for representations of bound quivers with monomial and commutativity relations over the virtual field $\mathbb{F}_1$, including the classification of indecomposable representations of a string pair over $\mathbb{F}_1$. In Section \ref{s:hall-algebra-f1}, we introduce the Hall algebra $\mathcal{H}_{\mathbb{F}_1}(Q,I)$ for a bound quiver $(Q,I)$ with monomial and commutativity relations. We prove that $\mathcal{H}_{\mathbb{F}_1}(Q,I)$ is a graded, connected, and cocommutative bialgebra (Proposition \ref{prop:hall-algebra-biaglebra}) and that its free subgroup generated by all indecomposable representations forms a Lie algebra (Corollary \ref{cor:hall-Lie-subalgebra-f1}). Section \ref{s:hall-algebra-nakayama-quiver} is devoted to comparing the Hall algebra $\mathcal{H}_{\mathbb{F}_1}(Q,I)$ with the degenerate Hall algebra $\mathcal{H}_1(Q,I)$ for a Nakayama bound quiver $(Q,I)$. In Theorem \ref{thm:equality-hall-numbers}, we prove that the equality in Question \ref{ques:tits-philosophy} (1) holds, and in Theorem \ref{thm:isom-hall-nakayma}, we establish the isomorphism between their Hall algebras. In Section \ref{s:hall-alg-gentle-one-cylce}, we address the class of gentle one-cycle bound quivers. In Theorem \ref{thm:hall-number-gentle-one}, we show that the equality in Question \ref{ques:tits-philosophy} (1) holds, and in Theorem \ref{thme:isom-hall-alg-gentle-one-cylce}, we establish the algebra isomorphism. We remark that our strategy for gentle one-cycle bound quivers differs from that for Nakayama bound quivers because Hall numbers are computed differently in the two settings.

\subsection*{Convention} Let $\mathbb{N}$ be the set of nonnegative integers and $\mathbb{N}^+$ the set of positive integers. For $a\geq b\in \mathbb{N}$, we denote by $\sqbinom{a}{b}$ the binomial coefficient and $\sqbinom{a}{b}_q$ the $q$-binomial coefficient.

% \subsection*{Use of AI Tools Declaration}
% The authors used Artificial Intelligence (AI) tools only to assist with language editing and the correction of grammatical errors. 

\section{Preliminaries}\label{s:preliminary}
\subsection{Hall algebras}
Let $\mathbb{F}_q$ be a finite field with $q$ elements and $\Lambda$ a finite-dimensional $\mathbb{F}_q$-algebra. Denote by $\Lambda$-${\rm mod}$ the category of finitely generated left $\Lambda$-modules, and $\opname{Iso}(\Lambda)$ a representative set of isomorphism classes of finitely generated left $\Lambda$-modules. For $M, L_1, \ldots, L_t\in \Lambda{\text{-{\rm mod}}}$, let $F_{L_1,\ldots,L_t}^M$ be the number of filtrations
\[
0=M_t\subset M_{t-1}\subset\cdots\subset M_1\subset M_0=M
\]
such that $M_{i-1}/M_i\cong L_i$ for $i=1,\ldots,t$. In particular, $F_{N,L}^M$ is the number of submodules $U$ of $M$ such that $U\cong L$ and $M/U\cong N$, which is called the {\em Hall number} of  $L,M,N$.

The following  formulas are very useful to compute Hall numbers.
\begin{lemma}\label{lem:hall-number-computation-k}
Let $L,M,N \in \Lambda\text{-{\rm mod}}$.

\begin{itemize}
    \item[(1)](Riedtmann--Peng formula \cite[Lemma 3.1]{P97}) \[
F^{M}_{N,L}
=
\frac{\left| \Ext^1_\Lambda(N,L)_M \right| \cdot a_M}
{\left| \Hom_\Lambda(N,L) \right| \cdot a_N \cdot a_L}.
\]Here $\Ext^1_\Lambda(N,L)_M$ 
in the subset  of $\Ext^1_\Lambda(N,L)$ consisting of 
equivalence classes of short exact sequences
$L \rightarrowtail M \twoheadrightarrow N$
with middle term isomorphic to $M$, and
$ a_X := \left| \operatorname{Aut}_\Lambda(X) \right|$ for any $X\in \Lambda\text{-{\rm mod}}$.
\item[(2)] Assume that $F_{N,L}^M\neq 0$, then
\[F_{N,L}^M=\frac{|\{f\in \Hom_{\Lambda}(L,M)\mid \operatorname{coker}f\cong N\}|}{a_L}=
  \frac{|\{g\in \Hom_{\Lambda}(M,N)\mid \operatorname{ker}g\cong L\}|}{a_N}.
   \]

\end{itemize}

\end{lemma}

The {\it Hall algebra} \cite{Ringel90banach} $\mathcal{H}(\Lambda)$ of $\Lambda$ is an associative algebra with unit whose underlying space is a free abelian group with a basis $\{u_{[M]}\mid M\in \opname{Iso}(\Lambda)\}$, and the multiplication is defined by
\[u_{[N]}u_{[L]}=\sum_{[M]}F_{N,L}^Mu_{[M]}.
\]
The following is also well known for experts, see, e.g., \cite[Lemma 5.2]{LW23}.
\begin{lemma}\label{lem:surjective-hall-algebra-homo}
    Let $A$ and $B$ be finite-dimensional $\mathbb{F}_q$-algebras such that there is a surjective homomorphism of algebras $\pi:A\to B$. Denote by $\iota: B\text{-$\opname{mod}$}\to A\text{-$\opname{mod}$}$ the embedding induced by $\pi$. For any triple $L,M,N\in B\text{-$\opname{mod}$}$, the equality $F_{N,L}^M=F_{\iota(N),\iota(L)}^{\iota(M)}$ holds. Consequently, there is a surjective homomorphism $f:\mathcal{H}(A)\to \mathcal{H}(B)$ given by \[u_{[M]}\mapsto\begin{cases}u_{[X]}&\text{if $M\cong \iota(X)$ for some $X\in B\text{-$\opname{mod}$}$};\\
    0&\text{else}.
    \end{cases}
    \]
\end{lemma}

\subsection{Hall polynomials}
Let $\mathbb{E}$ be a field extension of $\mathbb{F}_q$. For any $\mathbb{F}_q$-space $V$, we denote by $V^\mathbb{E}$ the $\mathbb{E}$-space $V\otimes_{\mathbb{F}_q}\mathbb{E}$.
Clearly, $\Lambda^\mathbb{E}$ is an $\mathbb{E}$-algebra, and for each $M\in\Lambda$-${\rm mod}$, $M^\mathbb{E}\in\Lambda^\mathbb{E}$-${\rm mod}$.  The field $\mathbb{E}$ is  {\it conservative} for $\Lambda$ if for any indecomposable $M\in \Lambda\text{-{\rm mod}}$, $(\End(M)/\rad \End(M))^\mathbb{E}$ is a field. Set
\[
\Omega(\Lambda)=\{\mathbb{E}\mid \text{$\mathbb{E}$ is a finite field extension of $\mathbb{F}_q$ which is conservative for $\Lambda$}\}.
\]
Now assume that $\Omega(\Lambda)$ is infinite  (e.g., $\Lambda$ is representation-finite). For $L,M,N\in\Lambda\text{-{\rm mod}}$, if there is a polynomial $\phi_{N,L}^M\in \mathbb{Z}[t]$ such that for any $\mathbb{E}\in \Omega(\Lambda)$, 
\[
\phi_{N,L}^M(|\mathbb{E}|)=F_{N^\mathbb{E},L^\mathbb{E}}^{M^\mathbb{E}},
\]
then $\phi_{N,L}^M(t)$ is called the {\it Hall polynomial}  associated with $L,M,N$. We say that $\Lambda$ has Hall polynomials, if any triple of $\Lambda$-modules $L,M,N$ has Hall polynomial $\phi_{N,L}^M(t)$.

The following is  a direct consequence of Lemma \ref{lem:surjective-hall-algebra-homo}.
\begin{corollary}\label{cor:hall-polynomial-factor-alg}
    Let $A$ and $B$ be finite-dimensional $\mathbb{F}_q$-algebras such that there is a surjective homomorphism of algebras $\pi:A\to B$. Assume that $\Omega(A)=\Omega(B)$. If $A$ has Hall polynomials, then $B$ also has Hall polynomials. Furthermore, for any triple $L,M,N\in B\text{-$\opname{mod}$}$, $\phi_{N,L}^M(t)=\phi_{\iota(N),\iota(L)}^{\iota(M)}(t)$.
\end{corollary}

It was conjectured by Ringel \cite{Ringel92} that every representation-finite $\mathbb{F}_q$-algebra has Hall polynomials. This conjecture has been verified for many important classes of representation-finite algebras. One approach in literature is to introduce certain partial order on the set of isomorphism classes of representations and prove the existence of Hall polynomials by induction with respect to the partial order. In particular, the following partial order by the length of extension groups is very useful, cf. \cite{GP97}.
\begin{lemma}\label{lem:order-relation-k}
  Let  \begin{tikzcd}[ampersand replacement=\&,cramped,sep=small]
L \arrow[r, tail, ] \& M \arrow[r, two heads, ] \& N 
\end{tikzcd} 
be a short exact sequence in $\Lambda{\text{-{\rm mod}}}$. Then 
\[
    \ell{(\Ext^1_\Lambda(M,M))}
    \leq \ell{(\Ext^1_\Lambda(L\oplus N,L\oplus N))},
\]
and equality holds if
and only if $M\cong L\oplus N $, where  $\ell{(\Ext^1_\Lambda(M,M))}$ is the length of $\Ext^1_\Lambda(M,M)$ as
$\mathbb{Z}$-module.
\end{lemma}

On the other hand, the Hall polynomials for some special triples of $\Lambda$-modules always exist, cf. \cite[Corollary 2.2.1]{P97}.
\begin{lemma}
    Let $\Lambda$ be a representation-finite $\mathbb{F}_q$-algebra and $S$ a simple $\Lambda$-module. Then, for any $M,L\in \Lambda{\text{-{\rm mod}}}$, the Hall polynomials $\phi_{S,L}^M(t)$ and $\phi_{L,S}^M(t)$ exist.
\end{lemma}
\begin{example}
    Let $Q$ be  the quiver $\xymatrix{1&2\ar[l]\ar[r]&3}$ and $A=\mathbb{F}_qQ$ the path algebra.  Denote by \[M=\xymatrix{\mathbb{F}_q&\mathbb{F}_q\ar[r]^1\ar[l]_1&\mathbb{F}_q},L_1=\xymatrix{\mathbb{F}_q&\mathbb{F}_q\ar[l]_1\ar[r]&0},\] 
    \[L_2=\xymatrix{0&\mathbb{F}_q\ar[l]\ar[r]^1&\mathbb{F}_q},S_2=\xymatrix{0&\mathbb{F}_q\ar[r]\ar[l]&0}.
    \] By using Lemma \ref{lem:hall-number-computation-k}, one can compute $\phi_{S_2,M}^{L_1\oplus L_2}(t)=t-1$.
\end{example}

\subsection{Degenerate Hall algebras and Lie subalgebras}

Now assume that $\Lambda$ admits Hall polynomials. 
The {\it generic Hall algebra}  $\mathcal{H}_t(\Lambda)$ of $\Lambda$ is an associative algebra whose underlying space is the free $\mathbb{Z}[t]$-module with basis $\{u_{[M]}\mid M\in \opname{Iso}(\Lambda)\}$, and the multiplication is given by 
\[u_{[N]}u_{[L]}=\sum_{[M]}\phi_{N,L}^M(t)u_{[M]}.
\]

Since $\phi_{N,L}^M(t)\in \mathbb{Z}[t]$, we may specialize $\phi_{N,L}^M(t)$ at $t=1$. The {\it degenerate Hall algebra} $\mathcal{H}_{1}(\Lambda)$ of $\Lambda$ is the associative algebra whose underlying space is the free abelian group with basis $\{u_{[M]}\mid M\in \opname{Iso}(\Lambda)\}$, and the multiplication is defined as
\[
u_{[N]}u_{[M]}=\sum_{[L]}\phi_{N,M}^L(1)u_{[L]}.
\]
Let $\operatorname{ind}\Lambda$ be a representative set of indecomposable $\Lambda$-modules.
Denote by $\mathfrak{n}_1(\Lambda)$  the free abelian subgroup of $\mathcal{H}_1(\Lambda)$ generated by $\{u_{[M]}\mid M\in \operatorname{ind}\Lambda\}$.
\begin{theorem}\cite{Ringel90banach,P97}
   The subgroup $\mathfrak{n}_1(\Lambda)$ is a Lie subalgebra of $\mathcal{H}_1(\Lambda)$, and $\mathcal{H}_1(\Lambda)\otimes_{\mathbb{Z}}\mathbb{Q}$ is the universal enveloping algebra of $\mathfrak{n}_1(\Lambda)\otimes_{\mathbb{Z}}\mathbb{Q}$.
\end{theorem}
We remark that the Lie bracket $[-,-]$ is given by \[[u_{[N]},u_{[L]}]:=\sum_{[M]\in \operatorname{ind}\Lambda}(\phi_{N,L}^M(1)-\phi_{L,N}^M(1))u_{[M]}.\]

\subsection{Representations of bound quivers over \texorpdfstring{$\mathbb{F}_1$}{F1}}

Following \cite{Sz12}, a finite-dimensional $\mathbb{F}_1$-vector space is a finite pointed set
$V=(V,0_V)$,
whose dimension is defined as
$\dim V=|V|-1$.
A map $f:V\to W$ between two finite-dimensional $\mathbb{F}_1$-vector spaces is called $\mathbb{F}_1$-linear if $f(0_V)=0_W$ and the restriction of $f$ to $V\setminus f^{-1}(0_W)$ is injective. The resulting category is denoted by $\operatorname{Vect} (\mathbb{F}_1)$.

Quiver representations over $\mathbb{F}_1$ were introduced in \cite{Sz12} and we refer to \cite{Sz12} for details. Obviously, quiver representations can be generalized to the setting of bound quivers with monomial and commutativity relations, cf. \cite{TF25}. 

Let $Q=(Q_0,Q_1,s,t)$ be a finite quiver, and let $I$ be a set of monomial and commutativity relations. We refer to such a pair $(Q,I)$ as a {\em bound quiver with monomial and commutativity relations}.

 An $\mathbb{F}_1$-representation
${V}=(V_i,V_\alpha)$
of $Q$ consists of an $\mathbb{F}_1$-vector space $V_i\in \operatorname{Vect} (\mathbb{F}_1)$ for each vertex $i\in Q_0$, and an $\mathbb{F}_1$-linear map
$V_\alpha:V_{s(\alpha)}\to V_{t(\alpha)}$
for each arrow $\alpha\in Q_1$. For any path $p=\alpha_k\cdots \alpha_2\alpha_1$ of length $k\geq 2$, we denote by 
\[
V_p=V_{\alpha_k}\circ\cdots\circ V_{\alpha_2}\circ V_{\alpha_1}:V_{s(\alpha_1)}\to V_{t(\alpha_k)}
\]
the composition of $V_{\alpha_1},\ldots, V_{\alpha_k}$.

 An $\mathbb{F}_1$-representation
of $(Q,I)$ is an $\mathbb{F}_1$-representation ${V}=(V_i,V_\alpha)$ of $Q$
such that
\begin{itemize}
    \item $V_p=0$ for each monomial relation $p\in I$;
    \item $V_p=V_q$ for each commutativity relation $p=q\in I$,
\end{itemize}
where $p,q$ are paths of length at least $2$. The dimension vector $\dimv {V}$ of ${V}$ is defined as
\[
\dimv {V}=(\dim V_i)_{i\in Q_0}\in \mathbb{Z}^{Q_0},
\]
and the dimension of ${V}$ is
\[\dim {V}:=\sum_{i\in Q_0}\dim V_i.\]

Given two $\mathbb{F}_1$-representations ${V}=(V_i,V_\alpha)$ and ${W}=(W_i,W_\alpha)$ of $(Q,I)$, a morphism $f=(f_i)_{i\in Q_0}:{V}\to {W}$ is a family of $\mathbb{F}_1$-linear maps $f_i:V_i\to W_i$ such that for every arrow $\alpha\in Q_1$, the following diagram commutes:
\[
\begin{tikzcd}
V_{s(\alpha)} \arrow[r,"V_\alpha"] \arrow[d,"f_{s(\alpha)}"'] 
& V_{t(\alpha)} \arrow[d,"f_{t(\alpha)}"] \\
W_{s(\alpha)} \arrow[r,"W_\alpha"'] 
& W_{t(\alpha)}.
\end{tikzcd}
\]
Moreover, if each $f_i$ is an isomorphism, then we say that $f$ is an isomorphism. Two representations ${V}$ and ${W}$ are isomorphic if there is an isomorphism $f$ from ${V}$ to ${W}$.
We denote by $\operatorname{rep}_{\mathbb{F}_1}(Q)$ the category of finite-dimensional $\mathbb{F}_1$-representations of $Q$, and  by $\operatorname{rep}_{\mathbb{F}_1}(Q,I)$
the full subcategory of $\operatorname{rep}_{\mathbb{F}_1}(Q)$ consisting of finite-dimensional $\mathbb{F}_1$-representations of $(Q,I)$. Denote by $\mathbf{0}$ the zero object of $\operatorname{rep}_{\mathbb{F}_1}(Q,I)$.

The
category $\operatorname{rep}_{\mathbb{F}_1}(Q,I)$ shares many good properties of the category $\operatorname{rep}_{\mathbb{F}_1}(Q)$. It has a zero object, kernels, cokernels, and direct sums. We refer to  \cite{Sz12} for the unexplained definitions. Moreover, the Jordan--H\"older theorem and the Krull--Schmidt theorem also hold for $\operatorname{rep}_{\mathbb{F}_1}(Q,I)$. The following is a direct consequence of the definition.
\begin{lemma}\label{lem:subrep-direct-sum-decomp}
    Let ${L}$ be a subrepresentation of ${M}=\bigoplus_{i=1}^t{M}_i$, where ${M}_i$ are direct summands of ${M}$. Then ${L}=\bigoplus_{i=1}^t({L}\cap {M}_i)$.
\end{lemma}

For any finite-dimensional $\mathbb{F}_1$-vector space $V=(V,0_V)$ and any field $\mathbb{F}$, one defines $V\otimes_{\mathbb{F}_1}\mathbb{F}$ to be the $\mathbb{F}$-vector space with basis $V\setminus\{0_V\}$. This gives a base change functor
\[
-\otimes_{\mathbb{F}_1}\mathbb{F}:\operatorname{Vect}(\mathbb{F}_1)\to \operatorname{Vect}(\mathbb{F}).
\]
Applying it vertexwise and arrowwise, one obtains a base change functor
\[
-\otimes_{\mathbb{F}_1}\mathbb{F}:\operatorname{rep}_{\mathbb{F}_1}(Q,I)\to \operatorname{rep}_\mathbb{F}(Q,I),
\]
where $\operatorname{rep}_\mathbb{F}(Q,I)$ is the category of finite-dimensional $\mathbb{F}$-representations of $(Q,I)$.
 This functor is exact and faithful, but in general not full. Moreover, it preserves direct sums.

\subsection{Representations of string pairs over $\mathbb{F}_1$}
In this subsection, we recall the basic construction of $\mathbb{F}_1$-representations for string pairs following \cite{TF25}.
\begin{definition}\cite{BR87}
Let $(Q,I)$ be a bound quiver with monomial relations.  We call $(Q,I)$  a {\em string pair} if $I$ is generated by monomial relations and the pair satisfies the following: 
\begin{itemize}
    \item[(s1)] for any $i \in Q_0$, there exist at most two arrows starting at $i$;
    \item[(s1')] for any $i \in Q_0$, there exist at most two arrows ending at $i$;
    \item[(s2)] for any $\alpha \in Q_1$, there exists at most one arrow $\beta \in Q_1$ such that  $\beta\alpha\notin I$;
    \item[(s2')] for any $\beta \in Q_1$, there exists at most one arrow $\alpha \in Q_1$ such that $\beta\alpha\notin I$.
\end{itemize}
\end{definition}

Let $(Q,I)$ be a string pair.
For an arrow $\beta\in Q_1$, we define a formal inverse $\beta^{-1}$ such that 
 $s(\beta^{-1})=t(\beta)$ and $t(\beta^{-1})=s(\beta)$. Denote by $Q_1^{-1}$ the set of formal inverses of arrows in $Q_1$.

A {\em walk} is a sequence $w=c_1\cdots c_n$, with $c_i\in Q_1\cup Q_1^{-1}$ for each $1\leq i\leq n$ and $s(c_i)=t(c_{i+1})$ for $1\leq i<n$. A {\em subwalk} of $w$ is a sequence of the form $c_i\cdots c_j$ for some $1\leq i\leq j\leq n$, and the inverse of $w$ is $c_n^{-1}\cdots c_{1}^{-1}$. We also denote by $s(w)=s(c_n)$ and $t(w)=t(c_1)$. A {\em string} of $(Q,I)$ of length $n$ is a walk $w=c_1\cdots c_n$ with $c_{i+1}\neq c_i^{-1}$ for $1\leq i<n$ and no subwalk of $w$ nor its inverse belong to $I$.
A {\em band} $b=b_1\dots b_n$ is a cyclic string, i.e., $s(b_n)=t(b_1)$, and such that any power $b^k$ is a string but $b$ itself is not a proper power of any string. For each $i\in Q_0$, we also introduce two (trivial) strings of length zero, $1_{i,1}$ and $1_{i,-1}$, whose sources and targets are $i$. Furthermore, $1_{i,\pm 1}^{-1}=1_{i,\mp 1}$.

Indecomposable representations of a string pair $(Q,I)$ over $\mathbb{F}_1$ were classified into string representations and band representations in \cite{TF25}. In contrast with the setting of fields, for a given band $b$ and a positive integer $m$, there exists a unique indecomposable band representation up to isomorphism.
In the sequel, since we will only consider string pairs which only admit string representations, we recall the construction and classification of string representations.

Let $u,v$ be strings of $(Q,I)$. We say $u\sim_{-1}v$ if $u=v$ or $u=v^{-1}$. It is clear that $\sim_{-1}$ is an equivalence relation on the set of strings of $(Q,I)$. Denote by $\mathbb{S}(Q,I)$ a set of representatives of $\sim_{-1}$ equivalence classes of strings.

\begin{definition}\label{def:string-rep}
Let $c=c_1\cdots c_n$ be a string of $(Q,I)$:
\[
\xymatrix{
u(0) & u(1)\ar[l]_{c_1} & \cdots\ar[l]_{c_2} & u(n)\ar[l]_{c_n}
},
\]
For each $i\in Q_0$, denote by 
\[
I_i:=\{k~|~u(k)=i, k=0,1,\dots, n\}.
\]
Let $\mathbf{V}_{\mathbb{F}_1}(c):=({V}_i,{V}_\alpha)$ be the $\mathbb{F}_1$-representation of $(Q,I)$ defined as follows:
\begin{itemize}
    \item ${V}_i$ is the $\mathbb{F}_1$-vector space of dimension $|I_i|$ with basis $\{z_k | k\in I_i\}$;
    \item  For any arrow $\alpha\in Q_1$, ${V}_\alpha:\{0,z_k~|~k\in I_{s(\alpha)}\}\to \{0,z_l~|~l\in I_{t(\alpha)}\}$ is given by
    \begin{align*}
    &{V}_\alpha(z_k)=z_{k-1} \quad \text{if $c_k=\alpha$};\\
    &{V}_\alpha(z_k)=z_{k+1} \quad \text{if $c_k=\alpha^{-1}$};\\
    &{V}_\alpha(z_k)=0\quad \text{else}.
    \end{align*}
\end{itemize}
The representation $\mathbf{V}_{\mathbb{F}_1}(c)$ is called  a {\em string representation} associated with $c$.
\end{definition}
\begin{remark}
Definition \ref{def:string-rep} is essentially the same as the definition of a string representation over a field $\mathbb{F}$. Namely, if one replaces $V_i$ by an $\mathbb{F}$-vector space of dimension $|I_i|$ with basis $\{z_k\mid k\in I_i\}$, then one obtains the string representation $\mathbf{V}_\mathbb{F}(c)$ of $(Q,I)$ over  $\mathbb{F}$. Clearly, $\mathbf{V}_\mathbb{F}(c)=\mathbf{V}_{\mathbb{F}_1}(c)\otimes_{\mathbb{F}_1}\mathbb{F}$.
\end{remark}

The following is a special case of \cite[Theorem 3.9]{TF25}, which is an analogue of a classical result for string pairs over fields \cite{BR87}.
\begin{proposition}\label{prop:classification-string-rep}
    For each string $c\in \mathbb{S}(Q,I)$, the string representation $\mathbf{V}_{\mathbb{F}_1}(c)$ is indecomposable. If $(Q,I)$ admits no band, then $\{\mathbf{V}_{\mathbb{F}_1}(c)\mid c\in \mathbb{S}(Q,I)\}$ forms a complete set of isomorphism classes of indecomposable representations of $(Q,I)$ over $\mathbb{F}_1$.
\end{proposition}

\section{Hall algebras for representations of bound quivers over $\mathbb{F}_1$}\label{s:hall-algebra-f1}

\subsection{Extension spaces}
Let $(Q,I)$ be a bound quiver with monomials and commutativity relations. Following Yoneda's construction, one defines the extension space $\Ext^i({N},{L})$ for any ${L},{N}\in \operatorname{rep}_{\mathbb{F}_1}(Q,I)$ and $i\geq 1$, cf. \cite{FRY24}. It is evident that $\Ext^i({N}, {L})$ possesses the structure of an $\mathbb{F}_1$-vector space. By definition, $\Ext^1({N}, {L})$ is finite-dimensional with the distinguished base point
\[
0=[
\begin{tikzcd}[ampersand replacement=\&, cramped, column sep=large]
{L} \arrow[r, tail, "{\sqbinom{1_{L}}{0}}"] \& {L}\oplus {N} \arrow[r, two heads, "{[0\ \ 1_{N}]}"] \& {N}
\end{tikzcd}
].
\]
In contrast, the higher extension spaces $\Ext^i({N}, {L})$ for $i \geq 2$ may be infinite-dimensional, see \cite{FYZ26a}.

  The following result was proved in \cite{FRY24} for linear quivers of type $\mathbb{A}$, and in \cite{GL26} for tree quivers. In fact, it holds for any bound quiver with monomial and commutativity relations.
\begin{lemma}\label{lem:blockdecomp}
Let  \begin{tikzcd}[ampersand replacement=\&,cramped,sep=small]
\epsilon: {L} \arrow[r, tail, "f"] \& {M} \arrow[r, two heads, "g"] \& {N}
\end{tikzcd} be a short exact sequence in $\operatorname{rep}_{\mathbb{F}_1}(Q,I)$. Suppose that ${M} = \bigoplus_{i=1}^s {M}_i$, where each ${M}_i$ is an indecomposable direct summand of ${M}$. Then $\epsilon$ is isomorphic to a direct sum of short exact sequences $\epsilon \cong \bigoplus_{i=1}^s \epsilon_i$, where each $\epsilon_i$ is of the form
\begin{tikzcd}[ampersand replacement=\&,cramped,sep=small]
 \widetilde{{L}}_i \arrow[r, tail, "f_i"] \& {M}_i \arrow[r, two heads, "g_i"] \& \widetilde{{N}}_i.
\end{tikzcd}
\end{lemma}

\begin{proof}
Let ${L}'$ be an arbitrary indecomposable direct summand of ${L}$.  By Lemma \ref{lem:subrep-direct-sum-decomp},
we have
\[
f({L}')=\bigoplus_{i=1}^s \bigl(f({L}')\cap {M}_i\bigr).
\]
Since $f$ is injective, $f({L}')\cong {L}'$ is indecomposable. It follows that there exists a unique index $j$ such that $f({L}')=f({L}')\cap {M}_j$.
For each $1\leq i\leq s$, let $\widetilde{{L}}_i$ be the maximal subrepresentation of ${L}$ such that $f(\widetilde{{L}}_i)\subseteq {M}_i$, and denote by $f_i=f|_{\widetilde{{L}}_i}$. Let $\epsilon_i$ be the short exact sequence \begin{tikzcd}[ampersand replacement=\&,cramped,sep=small]
 \widetilde{{L}}_i \arrow[r, tail, "f_i"] \& {M}_i \arrow[r, two heads, "g_i"] \& \operatorname{coker} f_i.
\end{tikzcd}
It follows that $\epsilon \cong \bigoplus_{i=1}^s \epsilon_i$.

\end{proof}

\begin{lemma}\label{lem:split}
Let
\begin{tikzcd}[ampersand replacement=\&,cramped,sep=small]
\epsilon: {L} \arrow[r, tail, "f"] \&{L}\oplus {N} \arrow[r, two heads, "g"] \& {N}
\end{tikzcd}
be a short exact sequence in $\operatorname{rep}_{\mathbb{F}_1}(Q,I)$,
then $\epsilon$ is split.
\end{lemma}

\begin{proof}
Let ${L}\oplus {N}=\bigoplus_{i=1}^s{M}_i$ be a decomposition of ${L}\oplus {N}$ into indecomposables. By Lemma \ref{lem:blockdecomp}, the short exact sequence $\epsilon$ is isomorphic to the direct sum of $\epsilon_1,\ldots, \epsilon_s$, where 
\begin{tikzcd}[ampersand replacement=\&,cramped,sep=small]
\epsilon_i: \widetilde{{L}}_i \arrow[r, tail, "f_i"] \& {M}_i\arrow[r, two heads, "g_i"] \& \widetilde{{N}}_i.
\end{tikzcd}
Since ${M}_i$ is nonzero, at least one of $\widetilde{{L}}_i$ and $\widetilde{{N}}_i$ is nonzero. On the other hand,  \[
\bigoplus_{i=1}^s\widetilde{L}_i\oplus \bigoplus_{i=1}^s\widetilde{N}_i\cong L\oplus N=\bigoplus_{i=1}^sM_i.
\]
It follows that, for each $i$, exactly one of $\widetilde{L}_i$ and $\widetilde{N}_i$ is nonzero, which turns out to be isomorphic to $M_i$.
Consequently, $\epsilon_i$ is split, and hence $\epsilon$ is split.

\end{proof}

\subsection{Hall algebra of $\operatorname{rep}_{\mathbb{F}_1}(Q,I)$}
In \cite{Sz12}, Szczesny introduced the Hall algebra of $\operatorname{rep}_{\mathbb{F}_1}(Q)$. Let $(Q,I)$ be a bound quiver with monomial and commutativity relations. 
In this subsection, we introduce the Hall algebra for $\operatorname{rep}_{\mathbb{F}_1}(Q,I)$ following \cite{Sz12}.

Since $\operatorname{rep}_{\mathbb{F}_1}(Q,I)$ is a  proto-exact category
in the sense \cite{DK19}, which is a non-additive generalization of Quillen’s exact category, one can define a Hall algebra $\mathcal{H}_{\mathbb{F}_1}(Q,I)$ for $\operatorname{rep}_{\mathbb{F}_1}(Q,I)$ in the same way as for $\Lambda$-${\rm mod}$.
Denote by $\operatorname{Iso}(Q,I)$ a representative set of isomorphism classes of objects in $\operatorname{rep}_{\mathbb{F}_1}(Q,I)$.

For any ${L},{M},{N}\in \operatorname{rep}_{\mathbb{F}_1}(Q,I)$, we denote by
$g_{{N},{L}}^{M}$ the number of subrepresentations ${U}$ of ${M}$ such that ${U}\cong {L}$ and ${M}/{U}\cong {N}$. 
\begin{lemma}\label{lem:hall-number-split-case}
    Let ${U}_1,\ldots, {U}_s\in \operatorname{rep}_{\mathbb{F}_1}(Q,I)$ be pairwise non-isomorphic indecomposable representations. Let ${L}=\bigoplus_{i=1}^sa_i{U}_i$, ${M}=\bigoplus_{i=1}^s(a_i+b_i){U}_i$ and ${N}=\bigoplus_{i=1}^sb_i{U}_i$, where $a_i,b_i$ are nonnegative integers. We have
    \[
    g_{{N},{L}}^{M}=\prod_{i=1}^s\sqbinom{a_i+b_i}{a_i}.
    \]
\end{lemma}
\begin{proof}
Considering ${L}$ as a subrepresentation of ${M}$, by  Lemma \ref{lem:subrep-direct-sum-decomp}, each indecomposable direct summand ${U}_i$ is a subrepresentation of exactly one indecomposable direct summand of ${M}$, say ${U}_j$. 
By Lemma \ref{lem:split}, the short exact sequence \begin{tikzcd}[ampersand replacement=\&,cramped,sep=small]
{L} \arrow[r, tail, "f"] \& {M} \arrow[r, two heads, "g"] \& {N}
\end{tikzcd} is split. It follows that $j=i$ and we conclude that $g_{{N},{L}}^{M}=\prod_{i=1}^s\sqbinom{a_i+b_i}{a_i}$.
\end{proof}

\begin{proposition}
    Let \[\mathcal{H}_{\mathbb{F}_1}(Q,I):=\bigoplus\limits_{{M}\in\operatorname{Iso}(Q,I) }\mathbb{Z}u_{[{M}]}\] be the free abelian group with a basis $\{u_{[{M}]}\mid {M}\in \opname{Iso}(Q,I)\}$. The following defines the structure of an associative algebra on $\mathcal{H}_{\mathbb{F}_1}(Q,I)$:
 \begin{equation}\label{eq:hall-multi}
u_{[{N}]}\star u_{[{L}]}=\sum_{[{M}]}g_{{N},{L}}^{{M}}u_{[{M}]},
\end{equation}
and the unit is given by $\eta:1\mapsto u_{[\mathbf{0}]}$.
\end{proposition}
\begin{proof}
    Note that there are finitely many representations of $(Q,I)$ up to isomorphism with a fixed dimension vector. It follows that the right hand side of \eqref{eq:hall-multi} is a finite sum, thus is well-defined.  The associativity follows the same proof as for Hall algebras over finite fields.
\end{proof}
We refer to $\mathcal{H}_{\mathbb{F}_1}(Q,I)$ as the {\it Hall algebra} of $\operatorname{rep}_{\mathbb{F}_1}(Q,I)$.
Let $K_0(\operatorname{rep}_{\mathbb{F}_1}(Q,I))$ be the Grothendieck group of $\operatorname{rep}_{\mathbb{F}_1}(Q,I)$. This yields a $K_0(\operatorname{rep}_{\mathbb{F}_1}(Q,I))$-grading on $\mathcal{H}_{\mathbb{F}_1}(Q,I)$. It is clear that the multiplication \eqref{eq:hall-multi} preserves the grading, hence $\mathcal{H}_{\mathbb{F}_1}(Q,I)$ becomes a $K_0(\operatorname{rep}_{\mathbb{F}_1}(Q,I))$-graded algebra.

For ${M}\in \operatorname{rep}_{\mathbb{F}_1}(Q,I)$ and $m\in \mathbb{N}^+$, denote 
\[
D_m({M})=\{({M}_1,\ldots,{M}_m)\mid {M}_1\oplus\cdots\oplus {M}_m\cong {M}\}\subset \operatorname{Iso}(Q,I)^m.
\]
Define a $\mathbb{Z}$-linear map $\Delta:\mathcal{H}_{\mathbb{F}_1}(Q,I)\to \mathcal{H}_{\mathbb{F}_1}(Q,I)\otimes \mathcal{H}_{\mathbb{F}_1}(Q,I)$ as follows
\[
\Delta(u_{[{M}]})=\sum_{({M}_1,{M}_2)\in D_2({M})}u_{[{M}_1]}\otimes u_{[{M_2}]},
\]
and $\epsilon:\mathcal{H}_{\mathbb{F}_1}(Q,I)\to \mathbb{Z}$ such that
\[
\epsilon(u_{[{M}]})=\begin{cases}
    0& \text{if ${M}\not\cong \mathbf{0}$};\\
    1& \text{else}.
\end{cases}
\]
\begin{proposition}\label{prop:hall-algebra-biaglebra}
    $(\mathcal{H}_{\mathbb{F}_1}(Q,I),\star, \Delta, \eta,\epsilon)$ is a graded, connected, and co-commutative bialgebra.
\end{proposition}
\begin{proof}
    The proof essentially follows the one in \cite{Sz12}.
    For completeness, we sketch the proof.   By the definition of $\Delta$, we know that $\Delta$ is compatible with the $K_0(\operatorname{rep}_{\mathbb{F}_1}(Q,I))$-grading and is co-commutative. The co-associativity of $\Delta$ follows from the fact that
    \[
    (\Delta\otimes \operatorname{id})\Delta(u_{[{M}]})=\sum_{({M}_1,{M}_2,{M}_3)\in D_3({M})}u_{[{M_1}]}\otimes u_{[{M}_2]}\otimes u_{[{M}_3]}=(\operatorname{id}\otimes \Delta)\Delta(u_{[{M}]}).
    \]
    For any ${M},{N}\in \operatorname{rep}_{\mathbb{F}_1}(Q,I)$,
    \begin{eqnarray}\label{eq:co-multi}
        \Delta(u_{[{M}]}\star u_{[{N}]})&=&\Delta(\sum_{\substack{({X},{R})\in Z_2(Q,I)\\ {X}\cong{N},{R}/{X}\cong{M}}}u_{[{R}]})\notag\\
        &=&\sum_{\substack{({X},{R})\in Z_2(Q,I)\\ {X}\cong{N},{R}/{X}\cong {M}}}(\sum_{({R}_1,{R}_2)\in D_2({R})}u_{[{R}_1]}\otimes u_{[{R}_2]}),
    \end{eqnarray}
    where $Z_2(Q,I):=\{({X},{R})\mid {R}\in \operatorname{Iso}(Q,I),{X}\subseteq {R}\}$, i.e., elements of $Z_2(Q,I)$ are pairs of an isomorphism class and a subobject.
    On the other hand, 
    \begin{eqnarray}\label{eq:multi-co}
        \Delta(u_{[{M}]})\star\Delta(u_{[{N}]})&=&(\sum_{({M}_1,{M}_2)\in D_2({M})}u_{[{M}_1]}\otimes u_{[{M}_2]})\star(\sum_{({N}_1,{N}_2)\in D_2({N})}u_{[{N}_1]}\otimes u_{[{N}_2]})\notag\\
        &=&\sum_{\substack{({M}_1,{M}_2)\in D_2({M})\\ ({N}_1,{N}_2)\in D_2({N})}}u_{[{M}_1]}\star u_{[{N}_1]}\otimes u_{[{M}_2]}\star u_{[{N}_2]}\notag\\
        &=&\sum_{\substack{({M}_1,{M}_2)\in D_2({M})\\ ({N}_1,{N}_2)\in D_2({N})}}\sum_{\substack{({X}_1,{R}_1)\in Z_2(Q,I)\\ {X}_1\cong {N}_1,{R}_1/{X}_1\cong {M}_1}}\sum_{\substack{({X}_2,{R}_2)\in Z_2(Q,I)\\ {X}_2\cong {N}_2,{R}_2/{X}_2\cong {M}_2}}u_{[{R}_1]}\otimes u_{[{R}_2]}.
    \end{eqnarray}
    For any $R_1$ and $R_2$, the coefficient of $u_{[R_1]}\otimes u_{[R_2]}$ in \eqref{eq:co-multi} is $g_{M,N}^{R_1\oplus R_2}$; while the coefficient of $u_{[R_1]}\otimes u_{[R_2]}$ in \eqref{eq:multi-co} is
    \[
    \sum_{\substack{M_1\oplus M_2\cong M\\ N_1\oplus N_2\cong N}}g_{M_1,N_1}^{R_1}g_{M_2,N_2}^{R_2}.
    \]
    
    By Lemma \ref{lem:subrep-direct-sum-decomp}, we obtain 
    \[
    g_{M,N}^{R_1\oplus R_2}=\sum_{\substack{M_1\oplus M_2\cong M\\ N_1\oplus N_2\cong N}}g_{M_1,N_1}^{R_1}g_{M_2,N_2}^{R_2}.
    \]
    We conclude that 
    \[
     \Delta(u_{[{M}]}\star u_{[{N}]})= \Delta(u_{[{M}]})\star\Delta(u_{[{N}]}).
    \]
    Finally, for any ${M}\in \operatorname{rep}_{\mathbb{F}_1}(Q,I)$,
    \begin{eqnarray*}
    (\epsilon\otimes \operatorname{id})\Delta(u_{[{M}]})&=&\epsilon\otimes \operatorname{id}(\sum_{({M}_1,{M}_2)\in D_2({M})} u_{[{M}_1]}\otimes u_{[{M}_2]})\\
    &=&\sum_{({M}_1,{M}_2)\in D_2({M})} \epsilon(u_{[{M}_1]}) u_{[{M}_2]}\\
    &=&\operatorname{id}(u_{[{M}]}).
    \end{eqnarray*}
    Similarly, one obtains $(\operatorname{id}\otimes \epsilon)\Delta= \operatorname{id}$.
    This finishes the proof.
\end{proof}
\begin{remark}
    There is a unique antipode $S$ on $\mathcal{H}_{\mathbb{F}_1}(Q,I)$ such that $\mathcal{H}_{\mathbb{F}_1}(Q,I)$ is a Hopf algebra over $\mathbb{Z}$.
\end{remark}

Let $\mathcal{H}(Q,I):=\mathcal{H}_{\mathbb{F}_1}(Q,I)\otimes_{\mathbb{Z}}\mathbb{C}$, and $\mathfrak{n}(Q,I)$  the subspace of $\mathcal{H}(Q,I)$ spanned by $\{u_{[{M}]}\mid \operatorname{ind}(Q,I)\}$, where $\operatorname{ind}(Q,I)$ is a representative set of indecomposable $\mathbb{F}_1$-representations of $(Q,I)$.
\begin{theorem}\label{thm:enveloping-lie-algebra}
    The subspace $\mathfrak{n}(Q,I)$ is a Lie subalgebra of $\mathcal{H}(Q,I)$, and $\mathcal{H}(Q,I)$ is isomorphic to the universal enveloping algebra of $\mathfrak{n}(Q,I)$.
\end{theorem}
\begin{proof}
   It is routine to show that $a\in \mathcal{H}(Q,I)$ is primitive if and only if $a\in \mathfrak{n}(Q,I)$. Consequently, $\mathfrak{n}(Q,I)$ is a Lie subalgebra of $\mathcal{H}(Q,I)$. The last assertion follows the Milnor--Moore theorem. 
\end{proof}
\begin{corollary}\label{cor:hall-Lie-subalgebra-f1}
    Let $\mathfrak{n}_{\mathbb{F}_1}(Q,I)$ be the free subgroup of $\mathcal{H}_{\mathbb{F}_1}(Q,I)$ generated by all the indecomposable representations. Then \[[u_{[{N}]},u_{[{L}]}]:=\sum_{{[M]}\in \operatorname{ind}(Q,I)}(g_{{N},{L}}^{M}-g_{{L},{N}}^{M})u_{[{M}]}\] defines a Lie algebra structure on $\mathfrak{n}_{\mathbb{F}_1}(Q,I)$.
\end{corollary}

\section{Hall algebras of Nakayama bound quivers}\label{s:hall-algebra-nakayama-quiver}

\subsection{Recollection on Nakayama algebras}\label{ss:structure-nakayama-algebra}
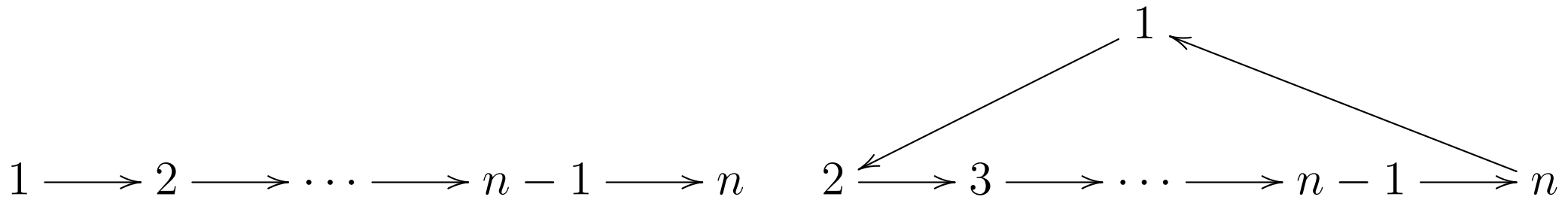
\begin{figure}[ht]
    \centering
    \begin{subfigure}{0.45\textwidth}
    \centering
        \xymatrix{1\ar[r]&2\ar[r]&\cdots \ar[r]&n-1\ar[r]&n}
        \caption{Linear quiver of type $\mathbb{A}_n$}
    \end{subfigure}
    \begin{subfigure}{0.45\textwidth}
        \xymatrix{&&1\ar[dll]&&\\
        2\ar[r]&3\ar[r]&\cdots\ar[r]&n-1\ar[r]&n\ar[ull]}
        \caption{Oriented cycle of length $n$}
    \end{subfigure}
   
    \caption{Ordinary quivers of Nakayama algebras}
    \label{fig:quiver-nakayama}
\end{figure}

Let $\mathbb{K}$ be an algebraically closed field and $A$ a finite-dimensional basic algebra over $\mathbb{K}$. The algebra $A$ is called a {\it Nakayama algebra} provided that $A$ is both right and left serial. Consequently,  every indecomposable $A$-module is uniserial. It is  well-known that every basic connected Nakayama algebra over $\mathbb{K}$ is a bound quiver algebra $\mathbb{K}Q_A/I$, where  $Q_A$ is either a linear oriented quiver of type $\mathbb{A}_n$ or an oriented cycle of length $n$ (cf. Figure \ref{fig:quiver-nakayama}), and $I$ is a set of monomial  relations. We will call the pair $(Q_A,I)$ a {\it linear/cyclic Nakayama bound quiver}. It is also known that the structure of representations of a Nakayama bound quiver is independent of the choice of the field.

Throughout this section, let $(Q,I)$ be a Nakayama bound quiver as in Figure \ref{fig:quiver-nakayama} and $\mathbb{F}_q$ a finite field with $q$ elements. It is evident that $(Q,I)$ is a string pair without bands. Consequently, the set $\{\mathbf{V}_{\mathbb{F}_q}(c) \mid c \in \mathbb{S}(Q,I)\}$ constitutes a complete set of isomorphism classes of indecomposable representations of $(Q,I)$ over $\mathbb{F}_q$. Similarly, according to Proposition \ref{prop:classification-string-rep}, the set $\{\mathbf{V}_{\mathbb{F}_1}(c) \mid c \in \mathbb{S}(Q,I)\}$ forms a complete set of isomorphism classes of indecomposable representations of $(Q,I)$ over $\mathbb{F}_1$.

For the study of Hall algebras, it is more convenient to adopt the following convention for indecomposable representations. Let $S_i$ denote the simple $\mathbb{F}_q$-representation associated with vertex $i$. Since every indecomposable representation of a Nakayama bound quiver over $\mathbb{F}_q$ is uniserial, it is uniquely determined by its top (or dually, its socle) and its length. We denote by $\mathbf{M}_{[i,l]}$ the unique indecomposable representation with top $S_i$ and length $l$ (whenever it exists); in particular, $S_i = \mathbf{M}_{[i,1]}$. For each vertex $i \in Q_0$, let $\ell_i$ be the maximal integer such that an indecomposable representation with top $S_i$ and length $\ell_i$ exists. 
 Then, the set$$\{\mathbf{M}_{[i,a_i]} \mid i \in Q_0, 1 \leq a_i \leq \ell_i\}$$forms a complete representative set of indecomposable representations of $(Q,I)$ over $\mathbb{F}_q$. For each $\mathbf{M}_{[i,a_i]}$, there exists a unique string $c \in \mathbb{S}(Q,I)$ such that $\mathbf{M}_{[i,a_i]} \cong \mathbf{V}_{\mathbb{F}_q}(c)$. To maintain uniformity across different base structures, we shall also denote the corresponding $\mathbb{F}_1$-representation as $\mathbf{M}_{[i,a_i]} := \mathbf{V}_{\mathbb{F}_1}(c)$ when no confusion arises.
The following well-known statement holds for both $\operatorname{rep}_{k}(Q,I)$ and $\operatorname{rep}_{\mathbb{F}_1}(Q,I)$.
\begin{lemma}\label{lem:canonical-sujective-to-simple}
    For any $i\in Q_0$ and $1\leq a\leq \ell_i$, $\dim \Hom(\mathbf{M}_{[i,a]},\mathbf{M}_{[i,1]})=1$ and the kernel of each nonzero homomorphism in $\Hom(\mathbf{M}_{[i,a]},\mathbf{M}_{[i,1]})$ is the radical of $\mathbf{M}_{[i,a]}$. Furthermore, if $i\neq j\in Q_0$, then $\dim\Hom(\mathbf{M}_{[i,a]},\mathbf{M}_{[j,1]})=0$.
\end{lemma}

\subsection{Hall polynomials of Nakayama algebras}

The existence of Hall polynomials for linear Nakayama algebras was shown by Ringel, since they are representation-directed. In the case of cyclic Nakayama algebras, the existence is a consequence of results by Guo \cite{G95} and Ringel \cite{Ringel93}, cf. \cite[Theorem 2.7]{G95} and \cite[Remark]{G95}.
\begin{theorem}\label{thm:exist-hall-nakayama}
    Let $(Q,I)$ be a cyclic Nakayama bound quiver and $A=\mathbb{F}_qQ/I$ the bound quiver algebra over $\mathbb{F}_q$. Then $A$ has Hall polynomials. 
\end{theorem}
\begin{proof}
Let us explain how Guo's work \cite{G95} implies the existence of Hall polynomials for cyclic Nakayama algebras. 
    Let $R$ be the arrow ideal of the path algebra $\mathbb{F}_qQ$ generated by all arrows. It has been shown by \cite[Theorem 2.7]{G95} that $\mathbb{F}_qQ/R^m$ has Hall polynomials for any $m\geq 0$. Now it is clear that there is a surjective homomorphism of algebra $f:\mathbb{F}_qQ/R^m\to A$ for some large enough integer $m$. Furthermore, $\Omega(A)=\Omega(\mathbb{F}_qQ/R^m)$. The result now follows directly from Corollary \ref{cor:hall-polynomial-factor-alg}.
\end{proof}

\begin{remark}\label{rem:cyclic-to-linear}
    In fact, we can also deduce the existence of Hall polynomials for linear Nakayama algebras from cyclic Nakayama algebras. Indeed, every linear Nakayama algebra with $n$ non-isomorphic simples is a quotient of a cyclic Nakayama algebra with $n+1$ non-isomorphic simples. The result now follows from Corollary \ref{cor:hall-polynomial-factor-alg}.
\end{remark}
\begin{remark}\label{rem:general-to-homogenous-m}
According to Remark \ref{rem:cyclic-to-linear} and the proof of Theorem \ref{thm:exist-hall-nakayama}, in order to calculate Hall numbers and Hall polynomials for a Nakayama bound quiver, it suffices to assume that $(Q,I)$ is a cyclic Nakayama bound quiver as in Figure \ref{fig:quiver-nakayama} (2) and $I=R^m$ for some positive integer $m\geq 2$. A similar reduction holds for the computation of Hall numbers of representations over $\mathbb{F}_1$.
\end{remark}

\begin{remark}
  In \cite{NI12}, Nasr-Isfahani announced an alternative proof of Theorem \ref{thm:exist-hall-nakayama}, but the argument appears to be incomplete. More precisely, it was claimed that every indecomposable module can be expressed as an integral linear combination of products of simple modules in $\mathcal{H}(A)$, and that the existence of Hall polynomials for $A$ follows from the existence of $\phi_{X,Y}^Z(t)$ when $X$ is simple. However, the correct statement is that every indecomposable module can be expressed as an integral linear combination of products of semisimple modules. Thus, to obtain the desired conclusion by this argument, one would need the existence of the Hall polynomials $\phi_{X,Y}^Z(t)$ for every semisimple module $X$.
\end{remark}
\begin{example}
   Let $Q$ be the quiver 
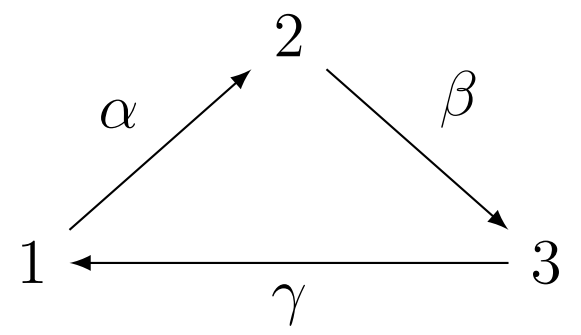
\[
\begin{tikzpicture}[>=latex]
\node (1) at (0,0) {$1$};
\node (2) at (1.7,1.5) {$2$};
\node (3) at (3.4,0) {$3$};

\draw[->] (1) -- node[above left] {$\alpha$} (2);
\draw[->] (2) -- node[above right] {$\beta$} (3);
\draw[->] (3) -- node[below] {$\gamma$} (1);
\end{tikzpicture}
\]
and $I=\{\beta\alpha\gamma, \alpha\gamma\beta\}$. We have
\[
    u_{[\mathbf{M}_{[3,3]}]}=u_{[S_3]}u_{[S_1]}u_{[S_2]}-u_{[S_1]}u_{[S_3]}u_{[S_2]}-u_{[S_3]}u_{[S_2]}u_{[S_1]}+u_{[S_1\oplus S_2 \oplus S_3]}.
\]
\end{example}

\begin{proposition}[{\cite[Proposition 3.4]{G95}}]\label{pro:semi-hall-number-k}
Let $(Q,I)$ be a cyclic Nakayama bound quiver as in Figure \ref{fig:quiver-nakayama}(2), where $I=R^m$ is the $m$-th power of arrow ideal $R$.
Let $M_1,M_2,M_3$ be representations of $(Q,I)$ over $\mathbb{F}_q$ such that $M_3$ is semisimple and $F_{M_3,M_1}^{M_2}\neq 0$. Assume that   $M_2=\bigoplus_{i=1}^{n} \bigoplus_{j=1}^{m}a_{i,j}^{(2)}\mathbf{M}_{[i,j]}$.
Then there exists a direct summand $N=\bigoplus_{i=1}^n\bigoplus_{j=1}^mb_{i,j}\mathbf{M}_{[i,j]}$ of $M_2$ such that $M_3\cong \operatorname{top} N$ and
\[
F_{M_3,M_1}^{M_2}
=\phi_{M_3,M_1}^{M_2}(q)=
\prod_{i,j}
\begin{bmatrix}
a_{i,j}^{(2)}\\
b_{i,j}
\end{bmatrix}_q
\,q^{\,b_{i,j}\left(\sum_{t=1}^{m-j} a_{i,j+t}^{(2)}-b_{i,j+t}\right)}.
\]

\end{proposition}

\subsection{The calculation of Hall numbers over $\mathbb{F}_1$}
In this section, we study the Hall numbers arising in the Hall algebra $\mathcal{H}_{\mathbb{F}_1}(Q,I)$ for  Nakayama bound quiver $(Q,I)$. Thanks to Remark \ref{rem:general-to-homogenous-m}, without loss of generality, we assume throughout this section that $(Q,I)$ is a cyclic Nakayama bound quiver as in Figure \ref{fig:quiver-nakayama}(2) with $I=R^m$ for some $m \ge 2$, where $R$ is the arrow ideal.

As we have explained in Section \ref{ss:structure-nakayama-algebra}, each representation ${M}$ in $\operatorname{rep}_{\mathbb{F}_1}(Q,I)$ and $\operatorname{rep}_{\mathbb{F}_q}(Q,I)$ can be written
uniquely as a direct sum of indecomposable representations
\[
{M}=\bigoplus_{i=1}^{n} \bigoplus_{j=1}^{m}a_{i,j}\mathbf{M}_{[i,j]}. 
\]
We will also write $\mathbf{a}=(a_{ij})\in \mathbb{Z}_{\geq 0}^{n\times m}$ for the representation ${M}$. Clearly, the length $\ell(\mathbf{a})$ of the representation $\mathbf{a}$ is $\sum_{i=1}^{n}\sum_{j=1}^mja_{i,j}$.
We define a total order $\prec_{1}$ on $\mathbb{N}^{n\times m}$ or, equivalently, on $\opname{Iso}(Q,I)$.

Assume that
$\mathbf{a}=(a_{i,j}),\mathbf{b}=(b_{i,j})\in \mathbb{N}^{n\times m}$.
Define $\mathbf{a}\prec_{1}\mathbf{b}$ if one of the following conditions is satisfied:
\begin{enumerate}
    \item
    $\ell(\mathbf{a})
    <
    \ell(\mathbf{b})$;
    \item
    $\ell(\mathbf{a})
    =
    \ell(\mathbf{b})$
    and there exists a pair $(s,t)$ such that
    $a_{s,t}<b_{s,t}$,and $ a_{i,j}=b_{i,j}$
    for all pairs $(i,j)$ satisfying
    $j>t$, or $j=t$  and  $ i<s$.
\end{enumerate}

\begin{lemma}\label{lem:semi-top-decompotion}
Let $M$ be  semisimple  and let $\psi:N\twoheadrightarrow M$ be an epimorphism in $\operatorname{rep}_{\mathbb{F}_1}(Q,I)$. Then there exists a direct sum decomposition $N=N_1\oplus N_2$ such that
$\operatorname{top} N_2\cong M$ and $\ker\psi=N_1\oplus \operatorname{rad} N_2$.
\end{lemma}
\begin{proof}
Let $K=\ker \psi$ and $N=X_1\oplus \cdots \oplus X_t$ be a decomposition of $N$ into indecomposable direct summands. Namely, each $X_j$ is isomorphic to $\mathbf{M}_{[i_j,l_j]}$ for some $i_j$ and $l_j$. By Lemma \ref{lem:blockdecomp}, the short exact sequence $K\rightarrowtail N\twoheadrightarrow M$ is isomorphic to the direct sum of short exact sequences $K_j\rightarrowtail X_j\twoheadrightarrow M_j$ for $1\leq j\leq t$. It follows that $M\cong M_1\oplus \cdots \oplus M_t$.
Since $M$ is semisimple,
we conclude that $M_j$ is either zero or isomorphic to $\mathbf{M}_{[i_j,1]}$ by Lemma \ref{lem:canonical-sujective-to-simple}. Without loss of generality, assume that $M_j$ is not zero for $1\leq j\leq s$. Let $N_2=X_1\oplus \cdots \oplus X_s$ and $N_1=X_{s+1}\oplus \cdots\oplus X_t$, then $\opname{top} N_2\cong M$ and $\ker \psi=N_1\oplus \opname{rad} N_2$.
\end{proof}

\begin{lemma}\label{lem:condition-hall-number-non-zero}
    Let $ \mathbf{a}^{(t)}=(a_{i,j}^{(t)})\in \mathbb{N}^{n\times m}$ for $t=1,2,3$, such that $\mathbf{a}^{(3)}$ is semisimple. Then $g_{\mathbf{a}^{(3)},\mathbf{a}^{(1)}}^{\mathbf{a}^{(2)}}\neq 0$ if and only if  there is a $\mathbf{b}=(b_{i,j})\in \mathbb{N}^{n\times m}$ satisfying the following conditions for all $i,j$:
    \begin{itemize}
        \item[(1)] $0\leq b_{i,j}\leq a_{i,j}^{(2)}$;
        \item[(2)] $\sum_{s=1}^mb_{i,s}=a_{i,1}^{(3)}$;
        \item[(3)] $a_{i,j}^{(2)}-a_{i,j}^{(1)}=b_{i,j}-b_{i-1,j+1}$.
    \end{itemize}
  Here, we adopt the convention that $b_{0,s} := b_{n,s}$, and $b_{i,j} = 0$ whenever $j > m$.
\end{lemma}
\begin{proof}
    
Assume $g_{\mathbf{a}^{(3)},\mathbf{a}^{(1)}}^{\mathbf{a}^{(2)}}\neq 0$. It follows that there exists an epimorphism $\psi:\mathbf{a}^{(2)}\twoheadrightarrow \mathbf{a}^{(3)}$. By Lemma \ref{lem:semi-top-decompotion}, there exist $\mathbf{b}=(b_{i,j})$ and $\mathbf{c}\in \mathbb{N}^{n\times m}$ such that $\mathbf{a}^{(2)}\cong \mathbf{b}\oplus \mathbf{c}$ and $\mathbf{a}^{(3)}\cong \mathbf{b}/\operatorname{rad} \mathbf{b}$. Here, $\mathbf{a}^{(2)}\cong \mathbf{b}\oplus \mathbf{c}$ implies condition (1), while $\mathbf{a}^{(3)}\cong \mathbf{b}/\operatorname{rad} \mathbf{b}$ implies condition (2). Condition (3) follows from the fact that $\operatorname{rad}\mathbf{b}=(b_{i-1,j+1})$.

Conversely, assume there exists a $\mathbf{b}=(b_{i,j})\in \mathbb{N}^{n\times m}$ satisfying conditions $(1)$–$(3)$. Let $c_{i,j}=a_{i,j}^{(2)}-b_{i,j}$ and set $\mathbf{c}:=(c_{i,j})$. By condition (1), $\mathbf{c}\in \mathbb{N}^{n\times m}$ and $\mathbf{a}^{(2)}\cong \mathbf{b}\oplus \mathbf{c}$. Condition (2) implies $\mathbf{a}^{(3)}\cong \mathbf{b}/\operatorname{rad}\mathbf{b}$. Thus, we have a canonical epimorphism $\psi:\mathbf{a}^{(2)}\cong \mathbf{b}\oplus \mathbf{c}\xrightarrow{[\pi,0]} \mathbf{b}/\operatorname{rad}\mathbf{b}\cong \mathbf{a}^{(3)}$. Finally, condition (3) implies that $\ker \psi\cong \mathbf{a}^{(1)}$, from which we conclude that $g_{\mathbf{a}^{(3)},\mathbf{a}^{(1)}}^{\mathbf{a}^{(2)}}\neq 0$. This completes the proof.

\end{proof}
\begin{remark}
\begin{itemize}
    \item[(1)]  The element  $\mathbf{b}$ in Lemma \ref{lem:condition-hall-number-non-zero} is uniquely determined by $\mathbf{a}^{(1)},\mathbf{a}^{(2)},\mathbf{a}^{(3)}$ via the conditions $(1)$-$(3)$. Moreover, $\mathbf{a}^{(3)}\cong \mathbf{b}/\operatorname{rad} \mathbf{b}$ and $\mathbf{a}^{(1)}\cong \operatorname{rad} \mathbf{b}\oplus (\mathbf{a}^{(2)}-\mathbf{b})$.
    \item[(2)]  Lemma \ref{lem:semi-top-decompotion} and Lemma \ref{lem:condition-hall-number-non-zero} also hold for the category $\operatorname{rep}_{\mathbb{F}_q}(Q,I)$, see \cite[Lemmas 3.1, 3.2]{G95}.
\end{itemize}
   
\end{remark}

\begin{proposition}
\label{pro:semi-hall-number-F_1}
Let $\mathbf{a}^{(t)}=(a_{i,j}^{(t)})\in \mathbb{N}^{n\times m}$ for $t=1,2,3$. Assume that $\mathbf{a}^{(3)}$ is semisimple and $g_{\mathbf{a}^{(3)},\mathbf{a}^{(1)}}^{\mathbf{a}^{(2)}}\neq 0$. Let $\mathbf{b}=(b_{i,j})\in \mathbb{N}^{n\times m}$ be determined by Lemma \ref{lem:condition-hall-number-non-zero}.
Then 
 \[
 g_{\mathbf{a}^{(3)},\mathbf{a}^{(1)}}^{\mathbf{a}^{(2)}}=\prod_{i,j}\sqbinom{a_{i,j}^{(2)}}{b_{i,j}}.
 \]
\end{proposition}

\begin{proof}
    Let $X$ be a subrepresentation of $\mathbf{a}^{(2)}$ such that $X\cong \mathbf{a}^{(1)}$ and $\mathbf{a}^{(2)}/X\cong \mathbf{a}^{(3)}$. According to Lemma \ref{lem:semi-top-decompotion}, there exist  subrepresentations $N_1$ and $N_2$ of $\mathbf{a}^{(2)}$  such that $\mathbf{a}^{(2)}=N_1\oplus N_2$ and the canonical projection $\pi:\mathbf{a}^{(2)}\twoheadrightarrow \mathbf{a}^{(2)}/X$ can be rewritten as
    \[
    N_1\oplus N_2\xrightarrow{[\pi_1,0]} N_1/\operatorname{rad} N_1=\mathbf{a}^{(2)}/X,
    \]
    where $\pi_1$ is the canonical projection, and $X=\operatorname{rad} N_1\oplus N_2$. According to the proof of Lemma \ref{lem:semi-top-decompotion},  $N_2$ is uniquely determined by $N_1$. We conclude that $X$ is in one-to-one correspondence with the choice of $N_1$. Note that $N_1\cong \mathbf{b}$ by the proof of Lemma \ref{lem:condition-hall-number-non-zero}. It follows that the number of choice of $N_1$ is $\prod_{i,j}\sqbinom{a_{i,j}^{(2)}}{b_{i,j}}$. This completes the proof.
\end{proof}

\begin{corollary}\label{cor:hall-number-semisimple}
    Let $\mathbf{a}^{(t)}=(a_{i,j}^{(t)})\in \mathbb{N}^{n\times m}$ for $t=1,2,3$, and assume that $\mathbf{a}^{(3)}$ is semisimple. Then $\phi_{\mathbf{a}^{(3)},\mathbf{a^{(1)}}}^{\mathbf{a}^{(2)}}(1)=g_{\mathbf{a}^{(3)},\mathbf{a^{(1)}}}^{\mathbf{a}^{(2)}}$.
\end{corollary}
\begin{proof}
    By \cite[Lemmas 3.1, 3.2]{G95}, we know that $\phi_{\mathbf{a}^{(3)},\mathbf{a^{(1)}}}^{\mathbf{a}^{(2)}}(t)\neq 0$ if and only if there exists a $\mathbf{b}$ satisfying the conditions $(1)$-$(3)$ in Lemma \ref{lem:condition-hall-number-non-zero}. Consequently, $\phi_{\mathbf{a}^{(3)},\mathbf{a^{(1)}}}^{\mathbf{a}^{(2)}}(t)=0$ if and only if $g_{\mathbf{a}^{(3)},\mathbf{a^{(1)}}}^{\mathbf{a}^{(2)}}=0$. The result now follows directly from Proposition \ref{pro:semi-hall-number-F_1} and Proposition \ref{pro:semi-hall-number-k}.
\end{proof}

In Theorem~\ref{thm:equality-hall-numbers}, we will show that the equality in Corollary~\ref{cor:hall-number-semisimple} actually holds for any triple of representations in $\operatorname{rep}_{\mathbb{F}_1}(Q,I)$. To this end, we first establish the following key lemma, which plays a crucial role in proving Theorem~\ref{thm:equality-hall-numbers}.

\begin{lemma}\label{lem:order-relation-nakayama-F_1}
    Let $\mathbf{a}=(a_{i,j})\in \mathbb{N}^{n\times m}$. Set $m_0=\max\{j\mid \exists\, a_{i,j}\neq 0\}$ and $\mathbf{b}=(b_{i,j})\in \mathbb{N}^{n\times m}$, where $b_{i,j}=\begin{cases}
        a_{i,m_0} &\text{if $j=m_0$};\\
        0&\text{else}.
    \end{cases}$
    Assume that there is a short exact sequence \[(\mathbf{a}-\mathbf{b})\oplus \operatorname{rad}\mathbf{b}\rightarrowtail \mathbf{a}'\twoheadrightarrow \operatorname{top}\mathbf{b}\] 
    in $\operatorname{rep}_{\mathbb{F}_1}(Q,I)$ or $\operatorname{rep}_{\mathbb{F}_q}(Q,I)$.
    If $\mathbf{a}'\neq \mathbf{a}$, then $\mathbf{a}'\prec_{1}\mathbf{a}$.
\end{lemma}
\begin{proof}

    Set $\mathbf{a}'=(a_{i,j}')\in \mathbb{N}^{n\times m}$.
We first note that $\ell(\mathbf{a})=\ell(\mathbf{a}')$. 

According to Lemma \ref{lem:semi-top-decompotion} (or \cite[Lemma 3.1]{G95} for $\operatorname{rep}_{\mathbb{F}_q}(Q,I)$), we may write $\mathbf{a}'=\mathbf{b}'\oplus \mathbf{c}'$ such that $\operatorname{top}\mathbf{b}'=\operatorname{top}\mathbf{b}$ and $\mathbf{c}'\oplus \operatorname{rad} \mathbf{b}'=(\mathbf{a}-\mathbf{b})\oplus \operatorname{rad} \mathbf{b}$, where $\mathbf{b}'=(b_{i,j}')$ and $\mathbf{c}'=(c_{i,j}')\in \mathbb{N}^{n\times m}$. By $\operatorname{top}\mathbf{b}'=\operatorname{top}\mathbf{b}$, we conclude that $\sum_{j=1}^mb_{i,j}'
=a_{i,m_0}$ for all $1\leq i\leq n$.
  A direct computation shows that 
\[
(\mathbf{a}-\mathbf{b})\oplus \operatorname{rad}\mathbf{b}=\begin{bmatrix}
    a_{1,1}&\cdots &a_{1,m_0-1}+a_{n,m_0}&0&\cdots &0\\
    \vdots &\cdots&\vdots&\vdots&\cdots &\vdots\\
    a_{i,1}&\cdots &a_{i,m_0-1}+a_{i-1,m_0}&0&\cdots&0\\
    \vdots &\cdots&\vdots&\vdots&\cdots &\vdots\\
    a_{n,1}&\cdots&a_{n,m_0-1}+a_{n-1,m_0}&0&\cdots&0
\end{bmatrix},
\]
\[
\mathbf{c}'\oplus \operatorname{rad}\mathbf{b}'=\begin{bmatrix}
    c_{1,1}'+b_{n,2}'&\cdots&c_{1,m-1}'+b_{n,m}'&c_{1,m}'\\
    \vdots&\cdots&\vdots&\vdots\\
    c_{i,1}'+b_{i-1,2}'&\cdots&c_{i,m-1}'+b_{i-1,m}'&c_{i,m}'\\
    \vdots&\cdots&\vdots&\vdots\\
    c_{n,1}'+b_{n-1,2}'&\cdots&c_{n,m-1}'+b_{n-1,m}'&c_{n,m}'
\end{bmatrix}.
\]
Note that all of $b_{ij}'$ and $c_{ij}'$ are nonnegative integers, we conclude that $c_{i,j}'=0$ for $j\geq m_0$ and $b_{i,l}'=0$ for $l\geq m_0+1$. In particular, the matrix $\mathbf{a}'$ has the following form
\[
\mathbf{a}'=\begin{bmatrix}
    \ast&\cdots &\ast& b_{1,m_0}'&0&\cdots &0\\
\vdots&\cdots&&\vdots&\cdots&\vdots\\
\ast&\cdots&\ast&b_{n,m_0}'&0&\cdots&0
\end{bmatrix}.
\]
For each $i$, $a_{i,m_0}=\sum_{j=1}^mb_{i,j}'\geq b_{i,m_0}'=a_{i,m_0}'$. If $a_{i,m_0}=a_{i,m_0}'$ for all $i=1,\ldots, n$, then $\mathbf{b}'=\mathbf{b}$. We conclude that $\mathbf{a}=\mathbf{a}'$ by 
 $(\mathbf{a}'-\mathbf{b})\oplus \operatorname{rad}\mathbf{b}=(\mathbf{a}-\mathbf{b})\oplus \operatorname{rad}\mathbf{b}$. 
Hence,  there exists an index $i$ such that $a_{i,m_0}'<a_{i,m_0}$, then $\mathbf{a}'\prec_1\mathbf{a}$ by definition. 
\end{proof}

We are now in a position to state the main result of this section.
\begin{theorem}\label{thm:equality-hall-numbers}
For any $\mathbf{a},\mathbf{b},\mathbf{c}\in \mathbb{N}^{n\times m}$,  we have $g_{\mathbf{c},\mathbf{a}}^{\mathbf{b}}=\phi_{\mathbf{c},\mathbf{a}}^{\mathbf{b}}(1)$.
\end{theorem}
\begin{proof}
We introduce a total order $\prec$ on $\mathbb{N}\times \mathbb{N}^{n\times m}$ as follows:
\[
(x,\mathbf{g})\prec (y,\mathbf{h}) \Leftrightarrow x<y \quad \text{or $x=y$ and $\mathbf{g}\prec_1\mathbf{h}$}.
\]
We will prove the statement by induction on $(\ell(\mathbf{b}), \mathbf{c})$. This is clearly true if $(\ell(\mathbf{b}), \mathbf{c})=(0,\mathbf{0})$. Assume that the statement is proved for any triple $\mathbf{a}',\mathbf{b}', \mathbf{c}'$ with $(\ell(\mathbf{b}'), \mathbf{c}')\prec (\ell(\mathbf{b}), \mathbf{c})$. Namely, $g_{\mathbf{c}',\mathbf{a}'}^{\mathbf{b}'}=\phi_{\mathbf{c}',\mathbf{a}'}^{\mathbf{b}'}(1)$
whenever $(\ell(\mathbf{b}'),\mathbf{c}')\prec (\ell(\mathbf{b}),\mathbf{c})$.

If $\mathbf{c}$ is semisimple, the result follows directly from Corollary \ref{cor:hall-number-semisimple}. Furthermore, if $\ell(\mathbf{b})\neq \ell(\mathbf{a})+\ell(\mathbf{c})$, then $g_{\mathbf{c},\mathbf{a}}^{\mathbf{b}}=0=\phi_{\mathbf{c},\mathbf{a}}^{\mathbf{b}}(1)$.
In the following, we assume that $\mathbf{c}$ is not semisimple and $\ell(\mathbf{b})=\ell(\mathbf{c})+\ell(\mathbf{a})$.

Let $m_0=\max\{j\mid \exists\, c_{ij}\neq 0\}$. Define $\mathbf{c}'=(c_{i,j}')$ and $\mathbf{c}''=(c_{i,j}'')$ as follows
\[
c_{i,j}'=\begin{cases}
    c_{i,j}&\text{if $j<m_0-1$};\\
    c_{i,m_0-1}+c_{i-1,m_0}&\text{if $j=m_0-1$};\\
    0&\text{if $j\ge m_0$}
\end{cases}\quad \text{and}\quad c_{i,j}''=\begin{cases}
    c_{i,m_0}&\text{if $j=1$};\\
    0&\text{else}.
\end{cases}
\]
By the construction, we have a short exact sequence $\mathbf{c}'\rightarrowtail \mathbf{c}\twoheadrightarrow \mathbf{c}''$. Consequently, $g_{\mathbf{c}'',\mathbf{c}'}^\mathbf{c}\neq 0$. Furthermore,  $g_{\mathbf{c}'',\mathbf{c}'}^\mathbf{c}=1$ by Lemma \ref{lem:condition-hall-number-non-zero} and Proposition \ref{pro:semi-hall-number-F_1}.

In the Hall algebra $\mathcal{H}_{\mathbb{F}_1}(Q,I)$, applying Lemma~\ref{lem:order-relation-nakayama-F_1} yields
\begin{align}\label{eq:hall-multi-expression}
    u_{[\mathbf{c}'']}\star u_{[\mathbf{c}']}&=\sum_{[\mathbf{d}]}g_{\mathbf{c}'',\mathbf{c}'}^\mathbf{d}u_{[\mathbf{d}]}\notag\\
    &=u_{[\mathbf{c}]}+\sum_{\substack{[\mathbf{d}] \\ \mathbf{d}\prec_1\mathbf{c}}} g_{\mathbf{c}'',\mathbf{c}'}^\mathbf{d} u_{[\mathbf{d}]}.
\end{align}
Multiplying \eqref{eq:hall-multi-expression} from the right by $u_{[\mathbf{a}]}$, we obtain
\begin{align}\label{eq:right-hand-hall-number}
    u_{[\mathbf{c}'']}\star u_{[\mathbf{c}']}\star u_{\mathbf{a}}&=u_{\mathbf{c}}\star u_{\mathbf{a}}+\sum_{\substack{[\mathbf{d}] \\ \mathbf{d}\prec_1\mathbf{c}}} g_{\mathbf{c}'',\mathbf{c}'}^\mathbf{d} u_{[\mathbf{d}]}\star u_{\mathbf{a}}\notag\\
    &=\sum_{[\mathbf{e}]}(g_{\mathbf{c},\mathbf{a}}^{\mathbf{e}}+\sum_{\substack{[\mathbf{d}] \\ \mathbf{d}\prec_1\mathbf{c}}} g_{\mathbf{c}'',\mathbf{c}'}^\mathbf{d}g_{\mathbf{d},\mathbf{a}}^{\mathbf{e}})u_{[\mathbf{e}]}.
\end{align}
On the other hand, associativity yields
\begin{align}\label{eq:left-hand-hall-number}
    u_{[\mathbf{c}'']}\star u_{[\mathbf{c}']}\star u_{\mathbf{a}}=\sum_{[\mathbf{e}]}(\sum_{[\mathbf{h}]}g_{\mathbf{c}'',\mathbf{h}}^{\mathbf{e}}g_{\mathbf{c}',\mathbf{a}}^{\mathbf{h}})u_{[\mathbf{e}]}.
\end{align}
Comparing the coefficients of $u_{[\mathbf{b}]}$ in \eqref{eq:right-hand-hall-number} and \eqref{eq:left-hand-hall-number}, we arrive at
\begin{align}\label{eq:recurse-in-F_1}
g_{\mathbf{c},\mathbf{a}}^{\mathbf{b}}=\sum_{\substack{[\mathbf{h}]\\ \ell(\mathbf{h})=\ell(\mathbf{c}')+\ell(\mathbf{a})}}g_{\mathbf{c}'',\mathbf{h}}^{\mathbf{b}}g_{\mathbf{c}',\mathbf{a}}^{\mathbf{h}}-\sum_{\substack{[\mathbf{d}] \\ \mathbf{d}\prec_1\mathbf{c}\\ \ell(\mathbf{d})=\ell(\mathbf{c})}} g_{\mathbf{c}'',\mathbf{c}'}^\mathbf{d}g_{\mathbf{d},\mathbf{a}}^{\mathbf{b}}.
\end{align}

By Corollary~\ref{cor:hall-number-semisimple}, we have $\phi_{\mathbf{c}'',\mathbf{c}'}^{\mathbf{c}}(1)=g_{\mathbf{c}'',\mathbf{c}'}^{\mathbf{c}}=1$. Thus, in the degenerate Hall algebra $\mathcal{H}_1(Q,I)$, we can expand the product as
\begin{align*}u_{[\mathbf{c}'']}\star u_{[\mathbf{c}']}&= \sum_{[\mathbf{d}]} \phi_{\mathbf{c}'',\mathbf{c}'}^\mathbf{d}(1) u_{[\mathbf{d}]}\\ &= u_{[\mathbf{c}]} + \sum_{\substack{[\mathbf{d}] \\ \mathbf{d}\prec_1\mathbf{c}}} \phi_{\mathbf{c}'',\mathbf{c}'}^\mathbf{d}(1) u_{[\mathbf{d}]}.
\end{align*}
Multiplying on the right by $u_{[\mathbf{a}]}$ and comparing the coefficients of $u_{[\mathbf{b}]}$, we obtain
\begin{align}\label{eq:recurse-in-k}\phi_{\mathbf{c},\mathbf{a}}^{\mathbf{b}}(1)= \sum_{\substack{[\mathbf{h}]\\ \ell(\mathbf{h})=\ell(\mathbf{c}')+\ell(\mathbf{a})}} \phi_{\mathbf{c}'',\mathbf{h}}^{\mathbf{b}}(1) \phi_{\mathbf{c}',\mathbf{a}}^{\mathbf{h}}(1) - \sum_{\substack{[\mathbf{d}] \\ \mathbf{d}\prec_1\mathbf{c}\\ \ell(\mathbf{d})=\ell(\mathbf{c})}} \phi_{\mathbf{c}'',\mathbf{c}'}^\mathbf{d}(1) \phi_{\mathbf{d},\mathbf{a}}^{\mathbf{b}}(1).
\end{align}
Since $\mathbf{c}$ is not semisimple and $\mathbf{c}''$ is a non-zero semisimple quotient of $\mathbf{c}$, it follows that $\ell(\mathbf{c}'') < \ell(\mathbf{c})$, and hence $(\ell(\mathbf{b}), \mathbf{c}'') \prec (\ell(\mathbf{b}), \mathbf{c})$. Moreover, the equality $\ell(\mathbf{h}) = \ell(\mathbf{c}') + \ell(\mathbf{a}) < \ell(\mathbf{b})$ implies $(\ell(\mathbf{h}), \mathbf{c}') \prec (\ell(\mathbf{b}), \mathbf{c})$. Similarly, we have $(\ell(\mathbf{d}), \mathbf{c}'') \prec (\ell(\mathbf{b}), \mathbf{c})$ and $(\ell(\mathbf{b}), \mathbf{d}) \prec (\ell(\mathbf{b}), \mathbf{c})$. In particular, each term appearing on the right-hand side of \eqref{eq:recurse-in-F_1} and \eqref{eq:recurse-in-k} strictly precedes $(\ell(\mathbf{b}), \mathbf{c})$ in the  order $\prec$. By the induction hypothesis, we have\[     g_{\mathbf{c}'',\mathbf{h}}^{\mathbf{b}} = \phi_{\mathbf{c}'',\mathbf{h}}^{\mathbf{b}}(1), \quad       g_{\mathbf{c}',\mathbf{a}}^{\mathbf{h}} = \phi_{\mathbf{c}',\mathbf{a}}^{\mathbf{h}}(1), \quad      g_{\mathbf{c}'',\mathbf{c}'}^{\mathbf{d}} = \phi_{\mathbf{c}'',\mathbf{c}'}^{\mathbf{d}}(1), \quad       g_{\mathbf{d}, \mathbf{a}}^{\mathbf{b}} = \phi_{\mathbf{d}, \mathbf{a}}^{\mathbf{b}}(1).\] Consequently, $g_{\mathbf{c}, \mathbf{a}}^{\mathbf{b}} = \phi_{\mathbf{c}, \mathbf{a}}^{\mathbf{b}}(1)$, completing the proof.

\end{proof}

\subsection{Isomorphism of Hall algebras}
As a consequence of Theorem~\ref{thm:equality-hall-numbers}, we obtain isomorphisms between the respective Hall algebras and their Lie subalgebras.
\begin{theorem}\label{thm:isom-hall-nakayma}
    Let $(Q,I)$ be a Nakayama bound quiver. The map $u_{[M]}\mapsto u_{[M]}$ defines an isomorphism of algebras from the degenerate Hall algebra $\mathcal{H}_1(Q,I)$ to the Hall algebra $\mathcal{H}_{\mathbb{F}_1}(Q,I)$. Furthermore, it restricts to an isomorphism of Lie algebras $\mathfrak{n}_{1}(Q,I)\cong \mathfrak{n}_{\mathbb{F}_1}(Q,I)$.
\end{theorem}
\begin{proof}
By the discussion in Section \ref{ss:structure-nakayama-algebra}, we know that $u_{[M]}\mapsto u_{[M]}$ yields an isomorphism of abelian groups. According to Remark \ref{rem:general-to-homogenous-m} and Theorem \ref{thm:equality-hall-numbers}, for any triple of representations $L,M,N$, we have $\phi_{N,L}^M(1)=g_{N,L}^M$. This equality of structure constants immediately implies that the map is an isomorphism of associative algebras. Finally, the Lie algebra isomorphism follows directly from the algebra isomorphism $\mathcal{H}_1(Q,I) \cong \mathcal{H}_{\mathbb{F}_1}(Q,I)$ and the definition of the Lie bracket.
    
\end{proof}

\section{Hall algebras of gentle one-cycle bound quivers}\label{s:hall-alg-gentle-one-cylce}
\subsection{Representations of gentle one-cycle algebras}
Let $\mathbb{F}_q$ be a finite field with $q$ elements. The gentle one-cycle algebra $\Lambda(n-1,1,1)$ is defined as a bound quiver algebra $\mathbb{F}_qQ/I$, where  $Q$ is a linearly oriented quiver of type $\mathbb{A}_n$ with  a loop at the vertex $n$ (cf. Figure \ref{fig:quiver-gentle one-cycle}), and $I=\{\alpha^2\}$. Throughout this section, we fix this bound quiver $(Q,I)$ and refer to it as a {\em gentle one-cycle bound quiver}.

\begin{figure}[ht]
  \centering
  $\xymatrix{
    1\ar[r] & 2\ar[r] & \cdots \ar[r] & n-1\ar[r] & n\ar@(ur,dr)^{\alpha}
  },\ \alpha^2$
  \caption{Bound quiver of the gentle one-cycle algebra}
  \label{fig:quiver-gentle one-cycle}
\end{figure}

% When $K$ is an algebraically closed field ,this is the derived-discrete algebra $\Lambda(n-1,1,1)$. We call the pair  $(Q_{BC},I)$ a {\it gentle one-cycle bound quiver}.

Clearly, $(Q,I)$ is a string pair and has no band. It follows from  \cite{BR87} and \cite{TF25} that both  the indecomposable representations in $\operatorname{rep}_{\mathbb{F}_q}(Q,I)$ and $\operatorname{rep}_{\mathbb{F}_1}(Q,I)$ are classified by strings. As a consequence, there is a one-to-one correspondence between isomorphism classes of   representations in $\operatorname{rep}_{\mathbb{F}_1}(Q,I)$ and isomorphism classes of representations in  $\operatorname{rep}_{\mathbb{F}_q}(Q,I)$.   If no confusion arises, we use the same notation $M$ to denote a representation $M$ in $\operatorname{rep}_{\mathbb{F}_1}(Q,I)$ with $M\otimes_{\mathbb{F}_1}\mathbb{F}_q$ in $\operatorname{rep}_{\mathbb{F}_q}(Q,I)$.
For instance, for any representations $M,N$, we write $\Ext^1_{\mathbb{F}_q}(M,N)$ for the first extension group of $M$ and $N$ as representations over $\mathbb{F}_q$ and $\Ext^1_{\mathbb{F}_1}(M,N)$ for the first extension space of $M$ and $N$ as representations over $\mathbb{F}_1$.

We recall the following explicit classifications of indecomposables in $\operatorname{rep}_{\mathbb{F}_q}(Q,I)$
given by \cite[Theorem 2.8]{CY24}. There are only finitely many indecomposable representations up to isomorphism.

\begin{theorem}[{\cite[Theorem 2.8]{CY24}}]\label{thm:indec-classification-gentle-one-cycle}
The following representations $U_{i,j}$, $V_i$ and $W_{i,j}$ form a complete set
of representatives of the isoclasses of indecomposable  representations of $(Q,I)$ over $\mathbb{F}_q$:

\begin{itemize}

\item $U_{i,j}$ for $1\le j\le i\le n$:
\[
{\small
\xymatrix@C=2.2em{
0 \ar[r] &
0 \ar[r] &
\cdots \ar[r] &
\overset{j}{\mathbb{F}_q} \ar[r]^{1} &
\cdots \ar[r]^{1} &
{\mathbb{F}_q} \ar[r]^{e_1} &
\overset{i}{\mathbb{F}_q^2} \ar[r]^{I} &
\cdots \ar[r]^{I} &
\overset{n}{\mathbb{F}_q^2} \ar@(ur,dr)^{M_\alpha}
}
}
\]

\item $U_{i,j}$ for $1\le i<j\le n$:
\[
{\small
\xymatrix@C=2.2em{
0 \ar[r] &
0 \ar[r] &
\cdots \ar[r] &
\overset{i}{\mathbb{F}_q} \ar[r]^{1} &
\cdots \ar[r]^{1} &
{\mathbb{F}_q} \ar[r]^{e_2} &
\overset{j}{\mathbb{F}_q^2} \ar[r]^{I} &
\cdots \ar[r]^{I} &
\overset{n}{\mathbb{F}_q^2} \ar@(ur,dr)^{M_\alpha}
}
}
\]

\item $V_i$ for $1\le i\le n$:
\[
{\small
\xymatrix@C=2.2em{
0 \ar[r] &
0 \ar[r] &
\cdots \ar[r] &
\overset{i}{\mathbb{F}_q} \ar[r]^{1} &
\cdots \ar[r]^{1} &
\overset{n}{\mathbb{F}_q} \ar@(ur,dr)^{0}
}
}
\]

\item $W_{i,j}$ for $1\le i\le j<n$:
\[
{\small
\xymatrix@C=2.2em{
0 \ar[r] &
0 \ar[r] &
\cdots \ar[r] &
\overset{i}{\mathbb{F}_q} \ar[r]^{1} &
\cdots \ar[r]^{1} &
\overset{j}{\mathbb{F}_q} \ar[r]^{0} &
0 \ar[r] &
\cdots \ar[r] &
\overset{n}{0} \ar@(ur,dr)^{0}
}
}
\]

\end{itemize}

Here
\[
e_1=
\begin{bmatrix}
1\\
0
\end{bmatrix},
\qquad
e_2=
\begin{bmatrix}
0\\
1
\end{bmatrix},
\qquad
I=
\begin{bmatrix}
1&0\\
0&1
\end{bmatrix},
\qquad
M_\alpha=
\begin{bmatrix}
0&0\\
1&0
\end{bmatrix}.
\]
\end{theorem}

In the following, we collect certain basic facts that will be frequently used in the subsequent sections.
\begin{lemma}\label{lem:basic-hom-space-gentle-one}
The following statements hold in $\operatorname{rep}_{\mathbb{F}_q}(Q,I)$ and $\operatorname{rep}_{\mathbb{F}_1}(Q,I)$:
    \begin{itemize}
        \item[(1)]  The simple representations are $S_1=W_{1,1},\ \ldots,\ S_{n-1}=W_{n-1,n-1}, S_n=V_n$;

\item[(2)]  For $W_{i,j}$, one has
\[
\operatorname{soc}(W_{i,j}) = S_j \quad \text{and}\quad \operatorname{top}(W_{i,j})=S_i.
\]
The nonzero indecomposable subrepresentations of $W_{i,j}$  are of the form $W_{r,j}$ $(i\le r\le j)$,
and the nonzero indecomposable quotient representations are of the form
$W_{i,l}$ $(i\le l\le j)$.

\item[(3)] For $V_i$, one has
\[
\operatorname{soc}(V_i) = S_n  \quad \text{and}\quad \operatorname{top} (V_i)=S_i.
\]
The nonzero indecomposable subrepresentations of $V_i$ are of the form $V_r$ $(i\leq r\leq n)$,
and the proper nonzero indecomposable quotient representations of $V_i$ are of the form
$W_{i,l}$ $(i\leq l <n)$.

\item[(4)]  For $U_{i,j}$, one has
\[
\operatorname{soc}(U_{i,j}) = S_n\quad \text{and}\quad \operatorname{top}(U_{i,j})=S_i\oplus S_j.
\]
In particular, $\dim \Hom(U_{i,i},S_i)=2$.
The nonzero indecomposable subrepresentations of $U_{i,j}$ are of type $U$ or $V$,
while the proper nonzero indecomposable quotient representations of $U_{i,j}$ are of type $V$ or $W$.

        \item[(5)] $\dim \Hom(U_{i,j},U_{i,j})=1$ if $j<i$ and $\dim\Hom(U_{i,j},U_{i,j})=2$ if $j\geq i$;
        \item[(6)] Let $\eta_{i1},\eta_{i2}\in \Hom(U_{i,i},S_i)$ be the homomorphism defined by $\begin{bmatrix}
            1\\0
        \end{bmatrix}$ and $\begin{bmatrix}
            0\\1
        \end{bmatrix}$, respectively. Let $\theta_i\in \Hom(U_{i,i},U_{i,i})$ be the homomorphism defined by $\begin{bmatrix}
            0&0\\1&0
        \end{bmatrix}$. Then $\eta_{i1}=\eta_{i2}\circ \theta_i$ and $\eta_{i1}\circ \theta_i=0$
    \end{itemize}
\end{lemma}
The existence of Hall polynomials for $\operatorname{rep}_{\mathbb{F}_q}(Q,I)$ has been recently established in \cite{CY24}.
\begin{theorem}\cite[Theorem 1.1]{CY24}
    Let $(Q,I)$ be a gentle one-cycle bound quiver and $\Lambda$ its bound quiver algebra over $\mathbb{F}_q$. Then $\Lambda$ has Hall polynomials.
\end{theorem}

\subsection{The calculation of Hall numbers over $\mathbb{F}_1$ }
In this subsection, we prove an analogue of Theorem \ref{thm:equality-hall-numbers} for gentle one-cycle bound quivers.
The inductive method used in the proof of Theorem \ref{thm:equality-hall-numbers} cannot be applied in this case, and we need to employ a new induction strategy, cf. \cite{GP97,Zhang}. We denote by $\Lambda=\mathbb{F}_qQ/I$ the bound quiver algebra over $\mathbb{F}_q$.
We begin with some special cases.

\begin{lemma}\label{lem:simple-hall-number-F1}
Let $L,M,S$ be representations of $(Q,I)$, where $S$ is simple. Then \[g_{S,L}^M=\phi_{S,L}^M(1).\]
\end{lemma}
\begin{proof}
We first compute the Hall number $F_{S,L}^M$ over $\mathbb{F}_q$. By Lemma \ref{lem:hall-number-computation-k} (2), 
\[
F^M_{S,L}
=
\frac{
\left|
\{\ f\in \operatorname{Hom}_\Lambda(M,S)\mid f\neq 0, \ker f\cong L\,\}
\right|
}
{q-1}.
\]
Without loss of generality, we may assume that $S=S_i$ and \[M=sU_{i,i}\oplus M_1\oplus \cdots \oplus M_t\oplus N,\]
where $s$ is a nonnegative integer,  $M_1,\ldots, M_t$ are indecomposable representations such that $\operatorname{top}(M_i)\cong S_i$, and $N$ is a representation such that $\operatorname{top}(N)$ does not have $S_i$ as a direct summand.  For each $1\leq j\leq t$, we fix a nonzero homomorphism $f_j\in \Hom(M_j,S_i)$, which forms a basis of $\Hom(M_j,S_i)$. By Lemma \ref{lem:basic-hom-space-gentle-one}, every morphism $f_\lambda$ in $\Hom(M,S)$ can be uniquely expressed as 
\[
f_{\lambda}=\begin{bmatrix}
    \lambda_{11}\eta_{i1}+\lambda_{12}\eta_{i2},\ldots, \lambda_{s1}\eta_{i1}+\lambda_{s2}\eta_{i2}, \lambda_1f_1,\ldots, \lambda_tf_t,0
\end{bmatrix},
\]
where $\lambda=(\lambda_{11},\lambda_{12},\ldots, \lambda_{s1},\lambda_{s2}, \lambda_1,\ldots, \lambda_t)\in \mathbb{F}_q^{2s+t}$. In particular, such a homomorphism $f_\lambda$ is uniquely determined by $\lambda\in \mathbb{F}_q^{2s+t}$. 

For $\lambda\in \mathbb{F}_q^{2s+t}$,
we write $\lambda(j)$ for the $j$-th component of $\lambda$ and define
\[
\operatorname{supp}(\lambda)=\{j\in \{1,\ldots, 2s+t\}\mid \lambda(j)\neq 0\}.
\]
For each $d\geq 1$, define
\[
c_d
=
\#\left\{
J\subseteq\{1,\dots,2s+t\}
\ \middle|\
|J|=d,\ \exists \ \lambda\in \mathbb{F}_q^{2s+t} \
\text{such that $\operatorname{supp}(\lambda)=J$ and $\ker f_\lambda\cong L$}
\right\}.
\]
We claim that if $\ker f_\lambda\cong L$ for some $\lambda\in \mathbb{F}_q^{2s+t}$, then for any $\lambda'\in \mathbb{F}_q^{2s+t}$ with $\operatorname{supp}(\lambda')=\operatorname{supp}(\lambda)$, $\ker f_{\lambda'}\cong L$. Indeed, let \[\varphi:sU_{i,i}\oplus M_1\oplus \cdots \oplus M_t\oplus N\xrightarrow{\operatorname{diag}\{
    \varphi_1, \cdots, \varphi_s, \varphi_{M_1},\cdots, \varphi_{M_t}, \varphi_N\}} sU_{i,i}\oplus M_1\oplus \cdots \oplus M_t\oplus N\]  be the diagonal homomorphism defined as follows:
\[
\varphi_N:=\mathbf{1}_N, \varphi_{M_i}:=\begin{cases}
    \mathbf{1}_{M_i}& \text{if $\lambda_i=0$};\\
    (\lambda_i')^{-1}\lambda_i\mathbf{1}_{M_i} &\text{if $\lambda_i\neq 0$}, 
\end{cases}
\]
\[
\varphi_l:=\begin{cases}
        \mathbf{1}_{U_{i,i}}&\text{if $\lambda_{l1}=\lambda_{l2}=0$};\\
        (\lambda_{l1}')^{-1}\lambda_{l1}\mathbf{1}_{U_{i,i}}&\text{if $\lambda_{l1}\neq 0$, $\lambda_{l2}=0$};\\
        (\lambda_{l2}')^{-1}\lambda_{l2}\mathbf{1}_{U_{i,i}}&\text{if $\lambda_{l1}= 0$, $\lambda_{l2}\neq 0$};\\
        (\lambda_{l2}')^{-1}\lambda_{l2}\mathbf{1}_{U_{i,i}}+(\lambda_{l2}')^{-1}(\lambda_{l1}-(\lambda_{l2}')^{-1}\lambda_{l2}\lambda_{l1}')\theta_i&\text{else}.
    \end{cases}
\]
Clearly, $\varphi$ is an isomorphism and the following diagram  commutes
\[
\xymatrix{M\ar[d]^{\varphi}\ar@{->>}[r]^{f_\lambda}&S\ar@{=}[d]\\
M\ar@{->>}[r]^{f_{\lambda'}}&S.}
\]
We conclude that $\ker f_\lambda\cong \ker f_{\lambda'}$.

According to the claim above, we have that
\[
F_{S,L}^M=\sum_{d=1}^{2s+t}c_d(q-1)^{d-1}.
\]
For any conservative finite extension $\mathbb{E}$ of $\mathbb{F}_q$, the same argument shows that $c_d$ is independent of the choice of $\mathbb{E}$ and
\[
F_{S^{\mathbb{E}},L^\mathbb{E}}^{M^{\mathbb{E}}}=\sum_{d=1}^{2s+t}c_d(|\mathbb{E}|-1)^{d-1}.
\]
It follows that $\phi_{S,L}^M(|\mathbb{E}|)\equiv c_1 \pmod{|\mathbb{E}|-1}$. As the cardinalities of the conservative finite extensions are unbounded,  we conclude that $\phi_{S,L}^M(1)=c_1$.

On the other hand,
\begin{align*}
    g_{S,L}^M&=|\{X\subseteq M\mid X\cong L, M/X\cong S\}|\\
    &=|\{f:M\twoheadrightarrow S\mid \ker f\cong L\}|.
\end{align*}
According to the definition of $\mathbb{F}_1$-linear map, each epimorphism $f:M\twoheadrightarrow S$ is uniquely determined by choosing a direct summand $S$ of $\operatorname{top}(M)$. It follows  that \[|\{f:M\twoheadrightarrow S\mid \ker f\cong L\}|=c_1,\] and hence $g_{S,L}^M=\phi_{S,L}^M(1)$.

\end{proof}

\begin{lemma}\label{lem:split-sum-hall-number-F1}
    For any representations $M,N$ of $(Q,I)$, $g_{M,N}^{M\oplus N}=\phi_{M,N}^{M\oplus N}(1)$.
\end{lemma}
\begin{proof}
We first compute the Hall number $F_{M,N}^{M\oplus N}$ over $\mathbb{F}_q$.

Let
\[
M=
\bigoplus_{1\leq i,j\leq n} a_{ij}U_{i,j}
\;\oplus\!
\bigoplus_{1\leq i\leq j<n} b_{ij}W_{i,j}
\;\oplus\!
\bigoplus_{1\leq i\leq n} c_iV_i,
\]
\[
N=
\bigoplus_{1\leq i,j\leq n} a'_{ij}U_{i,j}
\;\oplus\!
\bigoplus_{1\leq i\leq j<n} b'_{ij}W_{i,j}
\;\oplus\!
\bigoplus_{1\leq i\leq n} c'_iV_i.
\]
Since the middle term $M\oplus N$ corresponds to the split extension, one has
\[
\left|\Ext_{\mathbb{F}_q}^1(M,N)_{M\oplus N}\right|=1.
\]
Hence, by Riedtmann--Peng's formula, see Lemma \ref{lem:hall-number-computation-k},
\begin{align}\label{formula:ridt-peng-split}
    F_{M,N}^{M\oplus N}
=
\frac{ a_{M\oplus N}}
{\left| \Hom_\Lambda(M,N) \right| \cdot a_M \cdot a_N}.
\end{align}

For each indecomposable module $X\in \{U_{ij},W_{ij},V_i\}$, the ring
$\End_\Lambda(X)$ is local and
\[
\End_\Lambda(X)/\rad\End_\Lambda(X)\cong \mathbb{F}_q.
\]
Therefore, for every $m\geq 0$, the automorphism group of $mX$ is
the group of invertible $m\times m$ matrices over the local ring $\End_\Lambda(X)$, so
\[
|\operatorname{Aut}_\Lambda(mX)|
=
\prod_{s=1}^m
\Bigl(
|\End_\Lambda(X)|^m
-
|\End_\Lambda(X)|^{\,s-1}\,
|\rad\End_\Lambda(X)|^{\,m-(s-1)}
\Bigr).
\]
In our situation, we have
\[
|\End_\Lambda(U_{i,j})|=
\begin{cases}
q^2, & i\leq j,\\
q, & i> j,
\end{cases}
\qquad
|\rad\End_\Lambda(U_{i,j})|=
\begin{cases}
q, & i\leq j,\\
1, & i> j,
\end{cases}
\]
and
\[
|\End_\Lambda(W_{i,j})|=q,\qquad |\rad\End_\Lambda(W_{i,j})|=1,
\]
\[
|\End_\Lambda(V_i)|=q,\qquad |\rad\End_\Lambda(V_i)|=1.
\]
Hence, in each case, there exists an integer $\mu_X(m)\in\mathbb{N}$ such that
\[
|\operatorname{Aut}_\Lambda(mX)|
=
q^{\mu_X(m)}\prod_{s=1}^m (q^m-q^{s-1}).
\]
For convenience, write
\[
\gamma_m:=\prod_{s=1}^m (q^m-q^{s-1}).
\]
Then
\[
|\operatorname{Aut}_\Lambda(mX)|=q^{\mu_X(m)}\gamma_m.
\]
Note that morphisms between non-isomorphic indecomposable summands belong to radical, which contribute powers of $q$ in the calculation of automorphism groups. It follows that there exists an integer $\mu(M)\in \mathbb{N}$ such that
\[
|\operatorname{Aut}_\Lambda(M)|
=
q^{\mu(M)}
\prod_{1\leq i,j\leq n}\gamma_{a_{ij}}
\prod_{1\leq i\leq j<n}\gamma_{b_{ij}}
\prod_{1\leq i\leq n}\gamma_{c_i},
\]
Similarly,
\[
|\operatorname{Aut}_\Lambda(N)|
=
q^{\mu(N)}
\prod_{1\leq i,j\leq n}\gamma_{a'_{ij}}
\prod_{1\leq i\leq j<n}\gamma_{b'_{ij}}
\prod_{1\leq i\leq n}\gamma_{c'_i},
\]
and
\[
|\operatorname{Aut}_\Lambda(M\oplus N)|
=
q^{\mu(M\oplus N)}
\prod_{1\leq i,j\leq n}\gamma_{a_{ij}+a'_{ij}}
\prod_{1\leq i\leq j<n}\gamma_{b_{ij}+b'_{ij}}
\prod_{1\leq i\leq n}\gamma_{c_i+c'_i}.
\]
Moreover,
\[
|\Hom_\Lambda(M,N)|=q^{\dim_{\mathbb{F}_q}\Hom_\Lambda(M,N)}.
\]
Substituting these formulas into \eqref{formula:ridt-peng-split}, we get
\[
F_{M,N}^{M\oplus N}
=
q^{\nu(M,N)}
\prod_{1\leq i,j\leq n}
\frac{\gamma_{a_{ij}+a'_{ij}}}{\gamma_{a_{ij}}\gamma_{a'_{ij}}}
\prod_{1\leq i\leq j<n}
\frac{\gamma_{b_{ij}+b'_{ij}}}{\gamma_{b_{ij}}\gamma_{b'_{ij}}}
\prod_{1\leq i\leq n}
\frac{\gamma_{c_i+c'_i}}{\gamma_{c_i}\gamma_{c'_i}},
\]
for some integer $\nu(M,N)\in\mathbb{Z}$.
Since for $s,t\in \mathbb{N}$,
\[
\frac{\gamma_{s+t}}{\gamma_s\gamma_t}
=
q^{st}\sqbinom{s+t}{t}_q,
\]
it follows that
\[
F_{M,N}^{M\oplus N}
=
q^{\nu'(M,N)}
\prod_{1\leq i,j\leq n}\sqbinom{a_{ij}+a'_{ij}}{a'_{ij}}_q
\prod_{1\leq i\leq j<n}\sqbinom{b_{ij}+b'_{ij}}{b'_{ij}}_q
\prod_{1\leq i\leq n}\sqbinom{c_i+c'_i}{c'_i}_q
\]
for some 
$\nu'(M,N)\in\mathbb{Z}$.  Since each Gaussian binomial coefficient occurring above is congruent to $1$ modulo $q$, the product of these coefficients is not divisible by $q$. Since $F^{M\oplus N}_{M,N}$ is an integer, we must have $\nu'(M,N)\ge 0$.

By Lemma \ref{lem:hall-number-split-case},
\[
g_{M,N}^{M\oplus N}=\prod_{1\le i,j\le n}\sqbinom{a_{ij}+a'_{ij}}{a'_{ij}}
\prod_{1\le i\le j<n}\sqbinom{b_{ij}+b'_{ij}}{b'_{ij}}
\prod_{1\le i\le n}\sqbinom{c_i+c'_i}{c'_i}.
\]
It follows that  \[\phi_{M,N}^{M\oplus N}(q)\equiv g_{M,N}^{M\oplus N} \pmod{q-1}.\]
For any conservative finite extension $\mathbb{E}$ of $\mathbb{F}_q$, the same argument shows that
\[\phi_{M,N}^{M\oplus N}(|\mathbb{E}|)\equiv g_{M,N}^{M\oplus N} \pmod{|\mathbb{E}|-1}.
\]
As the cardinalities of the conservative finite extensions are unbounded, it follows that $\phi_{M,N}^{M\oplus N}(1)=g_{M,N}^{M\oplus N}$.

\end{proof}

For a representation $M\in \operatorname{rep}_{\mathbb{F}_q}(Q,I)$, we write $d(M)$ for the number of indecomposable direct summands of $M$, counted with multiplicities, and $\delta(M)=\ell(\Ext^1_{\mathbb{F}_q}(M,M))$ the length of $\Ext^1_{\mathbb{F}_q}(M,M)$ as $\mathbb{Z}$-module.

\begin{lemma}\label{lem:hall-number-equ when delta(N)=0}
    Assume that $\phi_{N,L}^M(1)=g_{N,L}^M$ for all $L,M,N\in \operatorname{rep}_{\mathbb{F}_q}(Q,I)$ with $N$ indecomposable. Then $\phi_{N,L}^M(1)=g_{N,L}^M$ for all $L,M,N\in \operatorname{rep}_{\mathbb{F}_q}(Q,I)$ with $\delta(N)=0$.
\end{lemma}
\begin{proof}
    We prove the statement by induction on $d(N)$. The statement is clearly true for $N=0$. If $d(N)=1$, then $N$ is indecomposable, and the assertion follows from the assumption.

    Now assume that $d(N)>1$ and write $N=N_1\oplus N_2$, where $N_1$ is indecomposable. Note that $\delta(N)=0$ implies that
    \begin{align*}
        0&=\Ext_{\mathbb{F}_q}^1(N_1\oplus N_2, N_1\oplus N_2)\\
        &=\Ext_{\mathbb{F}_q}^1(N_1,N_1)\oplus \Ext_{\mathbb{F}_q}^1(N_1,N_2)\oplus \Ext_{\mathbb{F}_q}^1(N_2, N_1)\oplus \Ext_{\mathbb{F}_q}^1(N_2,N_2).
    \end{align*}
    In particular, $\Ext_{\mathbb{F}_q}^1(N_1, N_2)=0$ and $\delta(N_1)=\delta(N_2)=0$.
     Note that $\Ext_{\mathbb{F}_q}^1(N_1,N_2)=0$ implies that $\Ext_{\mathbb{F}_1}^1(N_1,N_2)=0$. Indeed, if $\Ext_{\mathbb{F}_1}^1(N_1,N_2)\neq 0$, there is a short exact sequence 
     \begin{align}\label{seq:gentle-one-cylce}
         N_2\rightarrowtail X\twoheadrightarrow N_1
     \end{align} in $\operatorname{rep}_{\mathbb{F}_1}(Q,I)$ such that $X\not\cong N_1\oplus N_2$ by Lemma \ref{lem:split}. Note that the base change functor $-\otimes_{\mathbb{F}_1}\mathbb{F}_q$ is exact and induces a bijection between the isoclasses of representations in $\operatorname{rep}_{\mathbb{F}_1}(Q,I)$ and the isoclasses of representations in $\operatorname{rep}_{\mathbb{F}_q}(Q,I)$. The short exact sequence \eqref{seq:gentle-one-cylce} yields a non-split short exact sequence in $\operatorname{rep}_{\mathbb{F}_q}(Q,I)$, which contradicts  the fact that $\Ext^1_{\mathbb{F}_q}(N_1,N_2)=0$.
     
    By the associativity of both the Hall algebra $\mathcal{H}_{\mathbb{F}_1}(Q,I)$ and the degenerate Hall algebra $\mathcal{H}_1(Q,I)$, we obtain
    \begin{align}
        \sum_{[X]}g_{N_1,X}^Mg_{N_2,L}^X&=\sum_{[Y]}g_{N_1,N_2}^Yg_{Y,L}^M=g_{N_1,N_2}^{N}g_{N,L}^M,\label{eq:hall-induction-f1}\\
 \sum_{[X]}\phi_{N_1,X}^M(1)\phi_{N_2,L}^X(1)&=\sum_{[Y]}\phi_{N_1,N_2}^Y(1)\phi_{Y,L}^M(1)=\phi_{N_1,N_2}^{N}(1)\phi_{N,L}^M(1).\label{eq:hall-induction-k}
    \end{align}
  Note that $N_1$ is indecomposable, we have $g_{N_1,X}^M=\phi_{N_1,X}^M(1)$ and $g_{N_1,N_2}^N=\phi_{N_1,N_2}^N(1)$ by assumption. On the other hand, since $d(N_2)<d(N)$ and $\delta(N_2)=0$, we get $g_{N_2,L}^X=\phi_{N_2,L}^X(1)$ by the induction hypothesis. Furthermore, by Lemma \ref{lem:hall-number-split-case}, $g_{N_1,N_2}^N$ is nonzero. We conclude that  $g_{N,L}^M=\phi_{N,L}^M(1)$ by \eqref{eq:hall-induction-f1} and \eqref{eq:hall-induction-k}.
 
\end{proof}

\begin{proposition}\label{prop:equ when N is ind}
Assume that $g_{N,L}^M=\phi_{N,L}^M(1)$ for all $L,M,N\in \operatorname{rep}_{\mathbb{F}_q}(Q,I)$ with $N$ indecomposable. Then 
\[
\phi^{M}_{N,L}(1)=g^{M}_{N,L}
\]
holds for all $M,N,L\in \operatorname{rep}_{\mathbb{F}_q}(Q,I)$.
\end{proposition}
\begin{proof}
We denote by $<$ the lexicographic order on $\mathbb{N}\times \mathbb{N}^+$. Namely, $(l,d)<(l',d')$ if and only if $l<l'$ or $l=l'$ and $d<d'$. We prove  the statement by induction with respect to the order $<$.

If $d(N)=1$, the result holds by assumption. If $\delta(N)=0$, the statement follows from Lemma~\ref{lem:hall-number-equ when delta(N)=0}.

Assume that $\delta(N)>0$ and $d(N)>1$. Let $N=N_1\oplus N_2$ such that $N_1$ and $N_2$ are nonzero. It follows that $(\delta(N_1),d(N_1))<(\delta(N),d(N))$ and $(\delta(N_2),d(N_2))<(\delta(N),d(N))$. By the associativity of $\mathcal{H}_{\mathbb{F}_1}(Q,I)$, we obtain
\begin{align*}
    \sum_{[X]}g_{N_1,X}^Mg_{N_2,L}^X=\sum_{[Y]}g_{N_1,N_2}^Yg_{Y,L}^M=g_{N_1,N_2}^Ng_{N,L}^M+\sum_{\delta(Y)<\delta(N)}g_{N_1,N_2}^Yg_{Y,L}^M,
\end{align*}
where the second equality follows from Lemma \ref{lem:order-relation-k}. Equivalently,
\begin{align}\label{eq:hall-recursive-f11}
    g_{N_1,N_2}^Ng_{N,L}^M=\sum_{[X]}g_{N_1,X}^Mg_{N_2,L}^X-\sum_{\delta(Y)<\delta(N)}g_{N_1,N_2}^Yg_{Y,L}^M.
\end{align}
Similarly, in the degenerate Hall algebra $\mathcal{H}_1(Q,I)$, we have
\begin{align}\label{eq:hall-recursive-k1}
    \phi_{N_1,N_2}^N(1)\phi_{N,L}^M(1)=\sum_{[X]}\phi_{N_1,X}^M(1)\phi_{N_2,L}^X(1)-\sum_{\delta(Y)<\delta(N)}\phi_{N_1,N_2}^Y(1)\phi_{Y,L}^M(1).
\end{align}
Note that $(\delta(N_i),d(N_i))<(\delta(N),d(N))$ for $i=1,2$, and $(\delta(Y),d(Y))<(\delta(N),d(N))$. By the induction hypothesis, we have \[
g_{N_1,N_2}^N=\phi_{N_1,N_2}^N(1),\quad g_{N_1,X}^M=\phi_{N_1,X}^M(1),\quad g_{N_2,L}^X=\phi_{N_2,L}^X(1),\]
\[g_{N_1,N_2}^Y=\phi_{N_1,N_2}^Y(1),\quad g_{Y,L}^M=\phi_{Y,L}^M(1).
\]
By Lemma \ref{lem:hall-number-split-case}, we know that $g_{N_1,N_2}^N\neq 0$. Hence, we conclude that $g_{N,L}^M=\phi_{N,L}^M(1)$ from \eqref{eq:hall-recursive-f11} and \eqref{eq:hall-recursive-k1}.

\end{proof}

\begin{theorem}\label{thm:hall-number-gentle-one}
For all $L,M,N\in\operatorname{rep}_{\mathbb{F}_q}(Q,I)$, $\phi_{N,L}^{M}(1)=g_{N,L}^{M}$. 
\end{theorem}
\begin{proof}
    By Proposition \ref{prop:equ when N is ind} , it suffices to prove that for any indecomposable representation
$N$, $\phi_{N,L}^{M}(1)=g_{N,L}^{M}$ for any  $L,M\in \operatorname{rep}_{\mathbb{F}_q}(Q,I)$. We divide the proof into cases according to the type of $N$ in the classification of Theorem \ref{thm:indec-classification-gentle-one-cycle}.

\noindent{\bf Case 1: $N=W_{i,j},1\leq i\leq j <n$.} We prove the statement by induction on its length $\ell(N)$.  Note that $\ell(N)=1$ if and only if $i=j$. In this case, $N$
 is a simple representation, and the result holds directly by Lemma \ref{lem:simple-hall-number-F1}. 
Now assume that $j>i$. By Lemma \ref{lem:split-sum-hall-number-F1} and Lemma \ref{lem:hall-number-split-case}, we have
\[
g_{W_{i,i},W_{i+1,j}}^{W_{i,i}\oplus W_{i+1,j}}=\phi_{W_{i,i},W_{i+1,j}}^{W_{i,i}\oplus W_{i+1,j}}(1)=1\quad \text{and}\quad g_{W_{i+1,j},W_{i,i}}^{W_{i,i}\oplus W_{i+1,j}}=\phi_{W_{i+1,j},W_{i,i}}^{W_{i,i}\oplus W_{i+1,j}}(1)=1.
\]
Assume that there is a short exact sequence $W_{i+1,j}\rightarrowtail L\twoheadrightarrow W_{i,i}$ in $\operatorname{rep}_{\mathbb{F}_1}(Q,I)$ such that $L\not\cong W_{i,i}\oplus W_{i+1,j}$. According to Lemma \ref{lem:blockdecomp}, $L$ turns out to be indecomposable. By Theorem \ref{thm:indec-classification-gentle-one-cycle} and Lemma \ref{lem:basic-hom-space-gentle-one}, we conclude that $L\cong W_{i,j}$. Furthermore, $g_{W_{i,i},W_{i+1,j}}^{W_{i,j}}=\phi_{W_{i,i},W_{i+1,j}}^{W_{i,j}}(1)=1$ by Lemma \ref{lem:simple-hall-number-F1} and its proof. Similarly, one can show that if there is a short exact sequence $W_{i,i}\rightarrowtail L\twoheadrightarrow W_{i+1,j}$ in $\operatorname{rep}_{\mathbb{F}_1}(Q,I)$, then $L\cong W_{i,i}\oplus W_{i+1,j}$.

Putting all of these together, in the Hall algebra $\mathcal{H}_{\mathbb{F}_1}(Q,I)$, we have
\[
u_{[W_{i,i}]}u_{[W_{i+1,j}]}=u_{[W_{i,i}\oplus W_{i+1,j}]}+u_{[W_{i,j}]}, 
\]
\[
u_{[W_{i+1,j}]}u_{[W_{i,i}]}=u_{[W_{i,i}\oplus W_{i+1,j}]}.   
\]
Consequently,
\begin{align}\label{eq:hall-recur-wij}
u_{[W_{i,j}]}=u_{[W_{i,i}]}u_{[W_{i+1,j}]}-u_{[W_{i+1,j}]}u_{[W_{i,i}]}.
\end{align}
Multiplying \eqref{eq:hall-recur-wij} on the right by $u_{[L]}$ and comparing the coefficients of $u_{[M]}$, we obtain
\begin{align}\label{eq:rec-hall-wij-f1}
g^{M}_{W_{i,j},L}
=
\sum_{[X]}g^{M}_{W_{i,i},X}g^{X}_{W_{i+1,j},L}
-
\sum_{[Y]}g^{M}_{W_{i+1,j},Y}g^{Y}_{W_{i,i},L}.  
\end{align}
Similarly, one can show that $\eqref{eq:hall-recur-wij}$ also holds in the degenerate Hall algebra $\mathcal{H}_1(Q,I)$ (cf. the proof of \cite[Theorem 3.1]{CY24}). Hence one can obtain that 
\begin{align}\label{eq:rec-hall-wij-k}
\phi^{M}_{W_{i,j},L}(1)
=
\sum_{[X]}\phi^{M}_{W_{i,i},X}(1)\phi^{X}_{W_{i+1,j},L}(1)
-
\sum_{[Y]}\phi^{M}_{W_{i+1,j},Y}(1)\phi^{Y}_{W_{i,i},L}(1).  
\end{align}
By the induction hypothesis, we conclude that $g_{W_{i,j},L}^M=\phi_{W_{i,j},L}^M(1)$ from \eqref{eq:rec-hall-wij-f1} and \eqref{eq:rec-hall-wij-k}.

\noindent{\bf Case 2: $N=V_{i},1\leq i\leq n$.}
If  $i=n$, then $N=V_n$ is simple. The result holds by Lemma \ref{lem:simple-hall-number-F1}.
  For $1\leq i< n$, by the same argument as {\bf Case 1}, one can show that the following identities hold both in $\mathcal{H}_{\mathbb{F}_1}(Q,I)$ and $\mathcal{H}_1(Q,I)$:
\begin{align}
    u_{[W_{i,n-1}]}u_{[V_n]}&=u_{[W_{i,n-1}\oplus V_n]}+u_{[V_i]},\notag\\
     u_{[V_n]}u_{[W_{i,n-1}]}&=u_{[W_{i,n-1}\oplus V_n]},\notag\\
     u_{[V_i]}&=u_{[W_{i,n-1}]}u_{[V_n]}-u_{[V_n]}u_{[W_{i,n-1}]}.  \label{eq:hall-recur-vi}
\end{align}
Multiplying \eqref{eq:hall-recur-vi} on the right by $u_{[L]}$ and comparing the coefficients of $u_{[M]}$ in  $\mathcal{H}_{\mathbb{F}_1}(Q,I)$ and $\mathcal{H}_1(Q,I)$, we obtain
\begin{align*}
    g^{M}_{V_i,L}
&=
\sum_{[X]}g^{M}_{W_{i,n-1},X}g^{X}_{V_n,L}
-
\sum_{[Y]}g^{M}_{V_n,Y}g^{Y}_{W_{i,n-1},L},\\
\phi^{M}_{V_i,L}(1)
&=
\sum_{[X]}\phi^{M}_{W_{i,n-1},X}(1)\phi^{X}_{V_n,L}(1)
-
\sum_{[Y]}\phi^{M}_{V_n,Y}(1)\phi^{Y}_{W_{i,n-1},L}(1).
\end{align*}
We conclude $g_{V_i,L}^M=\phi_{V_i,L}^M(1)$ by {\bf Case 1} and Lemma \ref{lem:simple-hall-number-F1}.

\noindent{\bf Case 3: $N=U_{i,j},1\leq j \leq i \leq n$.}
Let us first consider the case $N=U_{n,j}$. Note that $U_{n,j}$ is projective in $\operatorname{rep}_{\mathbb{F}_q}(Q,I)$. It follows that $\Ext^1_{\mathbb{F}_q}(U_{n,j},-)=0$. Consequently, $\Ext^1_{\mathbb{F}_1}(U_{n,j},-)=0$ (cf. the proof of Lemma \ref{lem:hall-number-equ when delta(N)=0}). It follows that if $M\not\cong L\oplus U_{n,j}$, then $\phi_{U_{n,j},L}^M(1)=0=g_{U_{n,j},L}^M$. For the case $M\cong U_{n,j}\oplus L$, the equality holds by Lemma \ref{lem:split-sum-hall-number-F1}.

For $1\leq j\leq i< n$, by the same argument as {\bf Case 1}, one can show that the following identities hold both in $\mathcal{H}_{\mathbb{F}_1}(Q,I)$ and $\mathcal{H}_1(Q,I)$:
\begin{align*}
u_{[W_{i,n-1}]}u_{[U_{n,j}]}
&=
u_{[W_{i,n-1}\oplus U_{n,j}]}+u_{[U_{i,j}]},\\
u_{[U_{n,j}]}u_{[W_{i,n-1}]}
&=
u_{[W_{i,n-1}\oplus U_{n,j}]},\\
u_{[U_{i,j}]}
&=
u_{[W_{i,n-1}]}u_{[U_{n,j}]}-u_{[U_{n,j}]}u_{[W_{i,n-1}]}.
\end{align*}
Consequently, 
\begin{align*}
    g^{M}_{U_{i,j},L}
&=
\sum_{[X]}g^{M}_{W_{i,n-1},X}\,g^{X}_{U_{n,j},L}
-
\sum_{[Y]}g^{M}_{U_{n,j},Y}\,g^{Y}_{W_{i,n-1},L},\\
 \phi^{M}_{U_{i,j},L}(1)
&=
\sum_{[X]}\phi^{M}_{W_{i,n-1},X}(1)\,\phi^{X}_{U_{n,j},L}(1)
-
\sum_{[Y]}\phi^{M}_{U_{n,j},Y}(1)\,\phi^{Y}_{W_{i,n-1},L}(1),
\end{align*}
which implies that $g_{U_{i,j},L}^M=\phi_{U_{i,j},L}^M(1)$.

\noindent{\bf Case 4: $N=U_{i,j},1\leq i < j\leq n$.}
First, by the argument in {\bf Case 1}, one can show that the following identities hold both in $\mathcal{H}_{\mathbb{F}_1}(Q,I)$ and $\mathcal{H}_1(Q,I)$:
\begin{align*}
    u_{[V_i]}u_{[V_j]}&= u_{[V_i\oplus V_j]} + u_{[U_{j,i}]},\\
    u_{[V_j]}u_{[V_i]}&= u_{[V_i\oplus V_j]} + u_{[U_{i,j}]},\\
    u_{[U_{i,j}]}&= u_{[V_j]}u_{[V_i]}- u_{[V_i]}u_{[V_j]}+ u_{[U_{j,i}]}. 
\end{align*}
It follows that 
\begin{align*}
g^{M}_{U_{i,j},L}
=
\sum_{[X]}g^{M}_{V_j,X}g^{X}_{V_i,L}
-
\sum_{[Y]}g^{M}_{V_i,Y}g^{Y}_{V_j,L}
+
 g^{M}_{U_{j,i},L},\\
 \phi^{M}_{U_{i,j},L}(1)
=
\sum_{[X]}\phi^{M}_{V_j,X}(1)\phi^{X}_{V_i,L}(1)
-
\sum_{[Y]}\phi^{M}_{V_i,Y}(1)\phi^{Y}_{V_j,L}(1)
+
 \phi^{M}_{U_{j,i},L}(1).
\end{align*}
We conclude that $g_{U_{i,j},L}^M=\phi_{U_{i,j},L}^M(1)$ by {\bf Case 2} and {\bf Case 3}.

\end{proof}

\subsection{Isomorphism of Hall algebras}
As a consequence of Theorem \ref{thm:hall-number-gentle-one}, we obtain an isomorphism between the respective Hall algebras and their Lie subalgebras.
\begin{theorem}\label{thme:isom-hall-alg-gentle-one-cylce}
    Let $(Q,I)$ be a gentle one-cycle bound quiver. The map $u_{[M]}\mapsto u_{[M]}$ defines an isomorphism of algebras from the degenerate Hall algebra $\mathcal{H}_1(Q,I)$ to the Hall algebra $\mathcal{H}_{\mathbb{F}_1}(Q,I)$. Furthermore, it restricts to an isomorphism of Lie algebras $\mathfrak{n}_1(Q,I)\cong \mathfrak{n}_{\mathbb{F}_1}(Q,I)$.
\end{theorem}
\begin{remark}
    Following Riedtmann’s construction~\cite{Rie94}, one constructs a Lie algebra $\mathfrak{L}(Q,I)$ from $(Q,I)$. It was shown in~\cite{CY24} that $\mathfrak{L}(Q,I) \cong \mathfrak{n}_1(Q,I)$, which is further isomorphic to the positive part of the Lie algebra associated with a type $\mathbf{BC}$ root system (cf.~\cite{CY22}).
\end{remark}

\noindent{\bf Declaration of AI use}

The authors used ChatGPT solely to assist with grammar checking and language editing of the manuscript.

\bibliographystyle{plain}
\bibliography{myref}
\end{document}